\def\goth{\mathfrak}
\documentstyle[amssymb,12pt]{amsart}

\newtheorem{theorem}{Theorem}
\newtheorem{lemma}[theorem]{Lemma}

\newtheorem{definition}[theorem]{Definition}
\newtheorem{example}[theorem]{Examples}
\newtheorem{question}{Question}

\newcommand{\forces}{\Vdash}
\newcommand{\res}{\upharpoonright}
\newcommand{\ch}{$C\!H$}

\begin{document}

\baselineskip=18pt

  \begin{center}
     {\large Can a Small Forcing Create Kurepa Trees}
  \end{center}

  \begin{center}
 Renling Jin \& Saharon Shelah\footnote{The 
   	research of the second author was
	supported by  {\em Basic Research Foundation} of The Israel
	Academy of Sciences and Humanities.  Publication No 563.}

  \end{center}

  \bigskip

  \begin{quote}

    \centerline{Abstract}

    \small
In the paper we probe the possibilities of creating a Kurepa tree in a
generic extension of a model of \ch\, plus no Kurepa trees by an 
$\omega_1$-preserving forcing notion of size at most $\omega_1$. In
the first section we show that in the L\'{e}vy model obtained by
collapsing all cardinals between $\omega_1$ and a strongly inaccessible
cardinal by forcing with a countable support L\'{e}vy collapsing order
many $\omega_1$-preserving forcing notions of size at most $\omega_1$
including all $\omega$-proper forcing notions and some proper
but not $\omega$-proper forcing notions of size at most $\omega_1$ 
do not create Kurepa trees. In the second section we construct a model
of \ch\, plus no Kurepa trees, in which there is an $\omega$-distributive
Aronszajn tree such that forcing
with that Aronszajn tree does create a Kurepa tree in the generic
extension. At the end of the paper we ask three questions. 

\end{quote}

\section{Introduction}

By a model we mean a model of $Z\!F\!C$. By a forcing notion 
we mean a separative partially ordered set $\Bbb P$ with a largest 
element $1_{\Bbb P}$ used for a corresponding forcing extension.
Given a model $V$ of \ch, one can create a generic Kurepa tree by
forcing with an $\omega_1$-closed, $\omega_2$-c.c. forcing notion
no matter whether or not $V$ contains Kurepa trees [Je1]. One can
also create a generic Kurepa tree by forcing with a c.c.c. forcing
notion provided $V$ satisfies $\Box_{\omega_1}$ in addition [V].
Both forcing notions mentioned here have size at least $\omega_2$.
The size being at least $\omega_2$ seems necessary for guaranteeing
the generic trees have at least $\omega_2$ branches. On the other hand,
a Kurepa tree has a base set of size $\omega_1$, so it seems possible
to create a Kurepa tree by a forcing notion of size $\leqslant\omega_1$.
In this paper we discuss the following question: 
Given a model of \ch\, plus no Kurepa tree, whether can we 
find an $\omega_1$-preserving forcing notion of size $\leqslant\omega_1$ 
such that the forcing creates Kurepa trees?

This question is partially motivated by a parallel result about
Souslin tree. Given a ground model $V$. A Souslin tree
could be created by a c.c.c. forcing notion of size $\omega_1$ [ST]. 
There is also an $\omega_1$-closed forcing
notion of size $\omega_1$ which creates Souslin tree 
provided $V$ satisfies \ch\, [Je1]. The question
whether a Souslin trees could be created by a countable forcing notion
(equivalent to adding a Cohen real) turns out to be much harder. 
It was answered positively by the second author [S1] ten years ago.

We call a forcing notion $\omega_1$-preserving if $\omega_1$
in the ground model is still a cardinal in the generic extension.
In this paper we consider only $\omega_1$-preserving forcing notion
by the following reason. Let $V$ be the L\'{e}vy model. In $V$ there are no
Kurepa trees and \ch\, holds. Notice also that there is an 
$\omega_2$-Kurepa tree in $V$. If we simply collapse $\omega_1$ by forcing
with the collapsing order $Coll(\omega,\omega_1)$, the set of all
finite partial function from $\omega$ to $\omega_1$ ordered by reverse
inclusion, in $V$, then the $\omega_2$-Kurepa tree becomes a Kurepa tree
in $V^{Coll(\omega,\omega_1)}$. Notice also that $Coll(\omega,\omega_1)$
has size $\omega_1$ in $V$. So we require the forcing notions under
consideration be $\omega_1$-preserving to avoid the triviality.

In the first section we show some evidence that in the L\'{e}vy model
it is extremely hard to find a forcing notion, if it ever exists, of size
$\leqslant\omega_1$ which could create a Kurepa tree in the generic
extension. Assume our ground model $V$ is the  L\'{e}vy model. We show
first an easy result that any forcing notion of size 
$\leqslant\omega_1$ which adds no reals could not create Kurepa trees. 
Then we prove two main results: (1) For any stationary 
set $S\subseteq\omega_1$, 
if $\Bbb P$ is an $(S,\omega)$-proper forcing notion
of size $\leqslant\omega_1$, then there are no Kurepa trees in
the generic extension $V^{\Bbb P}$. Note that all axiom A forcing
notions are $(S,\omega)$-proper. (2) Some proper forcing notions
including the forcing notion for adding a club subset of $\omega_1$
by finite conditions do not create Kurepa trees in the generic
extension.

In the second section we show that there is a model of \ch\, plus no
Kurepa trees, in which there is an $\omega$-distributive Aronszajn
tree $T$ such that forcing
with $T$ does create a Kurepa tree in the generic extension.
We start with a model $V$ containing a strongly inaccessible cardinal 
$\kappa$. In $V$ we define an $\omega_1$-strategically closed, 
$\kappa$-c.c. forcing notion $\Bbb P$ such that forcing with 
$\Bbb P$ creates an $\omega$-distributive Aronszajn tree $T$ and a 
$T$-name $\dot{K}$ for a Kurepa tree $K$.
Forcing with $\Bbb P$ collapses also all cardinals between $\omega_1$
and $\kappa$ so that $\kappa$ is $\omega_2$ in $V^{\Bbb P}$.
Take $\bar{V} = V^{\Bbb P}$ as our ground model. 
Forcing with $T$ in $\bar{V}$ creates a Kurepa tree 
in the generic extension of $\bar{V}$. 
So the model $\bar{V}$ is what we are looking for except that
we have to prove that there are no Kurepa trees in $\bar{V}$, which 
is the hardest part of the second section.

We shall write $V$, $\bar{V}$, etc. for (countable) transitive 
models of $Z\!F\!C$.
For a forcing notion $\Bbb P$ in $V$ we shall write
$V^{\Bbb P}$ for the generic extension of $V$ by forcing with $\Bbb P$.
Sometimes, we write also $V[G]$ instead of $V^{\Bbb P}$
for a generic extension when a particular generic filter $G$
is involved. We shall fix a large enough regular cardinal
$\lambda$ throughout this paper and write $H(\lambda)$
for the collection of sets hereditarily of power less than 
$\lambda$ equipped with the membership relation. In a forcing argument
with a forcing notion $\Bbb P$ we shall write $\dot{a}$ for a
$\Bbb P$-name of $a$ and $\ddot{a}$ for a $\Bbb P$-name of
$\dot{a}$ which is again a $\Bbb Q$-name of $a$ for some forcing 
notion $\Bbb Q$. If $a$ is already in the ground model we
shall write simply $a$ for a canonical name of $a$. Let $\Bbb P$
be a forcing notion and $p\in {\Bbb P}$. We shall write
$q\leqslant p$ to mean $q\in {\Bbb P}$ and 
$q$ is a condition stronger than $p$. We shall
often write $p\forces$``$\ldots$'' for some $p\in {\Bbb P}$
instead of $p\forces_{\Bbb P}^{V}$``$\ldots$'' when the ground model $V$
and the forcing notion 
$\Bbb P$ in the argument is clear. We shall also write
$\forces$``$\ldots$'' instead of $1_{\Bbb P}\forces$``$\ldots$''. 
In this paper all of our trees are subtrees of the tree
$\langle 2^{<\omega_1},\subseteq\rangle$. So if $C$ is a linearly
ordered subset of a tree $T$, then $\bigcup C$ is the only possible
candidate of the least upper bound of $C$ in $T$.
In this paper all trees are growing upward. 
If a tree is used as a forcing notion we shall put the tree upside down. 
Let $T$ be a tree and $x\in T$. We write $ht(x)=\alpha$
if $x\in T\cap 2^{\alpha}$.  
We write $T_{\alpha}$ or $(T)_{\alpha}$, 
the $\alpha$-th level of $T$, for the set 
$T\cap 2^{\alpha}$ and write $T\!\res\!\alpha$ or $(T)\!\res\!\alpha$ for 
the set $\bigcup_{\beta<\alpha}T_{\beta}$.
We write $ht(T)$ for the height of $T$, which is the
smallest ordinal $\alpha$ such that $T_{\alpha}$ is empty.
By a normal tree we mean a tree $T$ such that (1) for any 
$\alpha<\beta<ht(T)$, for any $x\in T_{\alpha}$ there is an $y\in
T_{\beta}$ such that $x<y$; (2) for any $\alpha$ such that 
$\alpha +1 < ht(T)$ and for any $x\in T_{\alpha}$ there is $\beta<ht(T)$
and there are distinct $y_1, y_2\in T_{\beta}$ such that $x<y_1$ and
$x<y_2$. Given two trees $T$ and $T'$. We write $T\leqslant_{end}T'$
for $T'$ being an end-extension of $T$, {\em i.e.} $T'\res ht(T) = T$. 
By a branch of a tree $T$
we mean a totally ordered set of $T$ which intersects every
non-empty level of $T$. By an $\omega_1$-tree we mean a tree of height
$\omega_1$ with each of its levels at most countable.
A Kurepa tree is an $\omega_1$-tree with more than $\omega_1$ branches.
To see [J], [K] and [S2] for more information on
forcing, iterated forcing, proper forcing, etc. and 
to see [T] for more information on trees.

\bigskip

\noindent {\bf Acknowledgements}\quad The first part of this paper is
originated in 1993, when the first author was a Morrey assistant professor 
in University of California-Berkeley. He thanks H. Woodin 
for some inspiring discussion. The first author thanks also the
Department of Mathematics, Rutgers University for offering free housing
during his one week visit there at October, 1994, when the second part
of the paper was developed.

\section{Creating Kurepa Trees By a Small Forcing Is Hard}

First, we would like to state a theorem in [S2, 2.11] without proof
as a lemma which will be used in this section.
 
\begin{lemma}

In a model $V$ let $\Bbb P$ be a forcing notion and let $N$ be a 
countable elementary submodel of $H(\lambda)$. Suppose 
$G\subseteq {\Bbb P}$ is a $V$-generic filter. Then
\[N[G]=\{\dot{a}_G:\dot{a}\mbox{ is a }{\Bbb P}\mbox{-name and }\dot{a}\in N\}\]
is a countable elementary submodel of $(H(\lambda))^{V[G]}$.

\end{lemma}

We choose the L\'{e}vy model $\bar{V}=V^{Lv(\kappa,\omega_1)}$ 
as our ground model throughout
this section, where $\kappa$ is a strongly inaccessible cardinal
in $V$ and $Lv(\kappa,\omega_1)$, the Levy collapsing order,
is the set
\[\begin{array}{ll}
\{p\subseteq (\kappa\times\omega_1)\times\kappa:&p\mbox{ is 
a countable function and }\\
 & (\forall (\alpha,\beta)\in dom(p)) (p(\alpha,\beta)\in\alpha)\}
\end{array}\]
ordered by reverse inclusion. For any 
$A\subseteq\kappa$ we write $Lv(A,\omega)$ for the set of all
$p\in Lv(\kappa,\omega_1)$ such that $dom(p)\subset A\times\omega_1$.

We now prove an easy result.

\begin{theorem}

Let $\Bbb P$ be a forcing notion
of size $\leqslant\omega_1$ in $\bar{V}$. If forcing with $\Bbb P$ does
not add new countable sequences of ordinals, 
then there are no Kurepa trees in $\bar{V}^{\Bbb P}$.

\end{theorem}

\noindent {\bf Proof:}\quad
Since $\Bbb P$ has size $\leqslant\omega_1$, there is an $\eta<\kappa$
such that ${\Bbb P}\in V^{Lv(\eta,\omega_1)}$. Hence $\bar{V}^{\Bbb P}
=V^{(Lv(\eta,\omega_1)*\dot{\Bbb P})\times 
Lv(\kappa\smallsetminus\eta,\omega_1)}$.
But $Lv(\kappa\smallsetminus\eta,\omega_1)$ in $V$ is again a Levy
collapsing order in $V^{Lv(\eta,\omega_1)*\dot{\Bbb P}}$ because 
$\Bbb P$ adds no new countable sequences of ordinals, so that the forcing 
notion $Lv(\kappa\smallsetminus\eta,\omega_1)$ is absolute between $V$
and $V^{Lv(\eta,\omega_1)*\dot{\Bbb P}}$.
Hence there is no
Kurepa trees in $\bar{V}^{\Bbb P}$. \quad $\Box$

\medskip

Next we prove the results about $(S,\omega)$-proper forcing notions.

\begin{definition}

A forcing notion $\Bbb P$ is said to satisfies property ($\dagger$) if
for any $x\in H(\lambda)$, there exists a sequence 
$\langle N_i:i\in\omega\rangle$ of elementary
submodels of $H(\lambda)$ such that 

(1) $N_i\in N_{i+1}$ for every $i\in\omega$,

(2) $\{{\Bbb P}, x\}\subseteq N_0$,

(3) for every $p\in {\Bbb P}\cap N_0$ 
there exists a $q\leqslant p$ and $q$ is
$({\Bbb P}, N_i)$-generic for every $i\in\omega$.

\end{definition}
 
\begin{lemma}

Let $V$ be any model. Let $\Bbb P$ and $\Bbb Q$ be two forcing notions
in $V$ such that $\Bbb P$ has size $\leqslant\omega_1$ and satisfies property 
($\dagger$), and $\Bbb Q$ is $\omega_1$-closed (in $V$). 
Suppose $T$ is an $\omega_1$-tree in $V^{\Bbb P}$.
Then $T$ has no branches which are in $V^{{\Bbb P}\times {\Bbb Q}}$ but
not in $V^{\Bbb P}$.

\end{lemma}

\noindent {\bf Proof:}\quad
Suppose, towards a contradiction, that there is a branch $b$ of $T$ in 
$V^{{\Bbb P}\times {\Bbb Q}}\smallsetminus V^{\Bbb P}$.
Without loss of generality, we can assume that
\[\forces_{\Bbb P} \forces_{\Bbb Q} 
(\ddot{b}\mbox{ is a branch of }\dot{T}\mbox{ in }
V^{{\Bbb P}\times {\Bbb Q}}\smallsetminus V^{\Bbb P}).\]

\medskip

{\bf Claim 4.1}\quad For any $p\in {\Bbb P}$, $q\in {\Bbb Q}$,
$n\in\omega$ and 
$\alpha\in\omega_1$, there are $p'\leqslant p$, $q_j\leqslant q$  for
$j<n$ and $\beta\in\omega_1\smallsetminus\alpha$ such that 
\[ p'\forces((\exists \{t_j:j<n\}\subseteq\dot{T}_{\beta}) 
((j\not=j'\rightarrow t_j\not=t_{j'})\wedge
\bigwedge_{j<n}(q_j\forces t_j\in\ddot{b}))).\]

Proof of Claim 4.1:\quad
Since 
\[p\forces_{\Bbb P} q\forces_{\Bbb Q} 
(\ddot{b}\mbox{ is a branch of }\dot{T}\mbox{ in }
V^{{\Bbb P}\times {\Bbb Q}}\smallsetminus V^{\Bbb P}),\]
then $p$ forces that $q$ can't determine $\ddot{b}$. Hence
\[ p\forces((\exists\beta\in\omega_1\smallsetminus\alpha)
(\exists q_j\leqslant q\mbox{ for }j<n) (\exists t_j\in\dot{T}_{\beta} 
\mbox{ for }j<n )\]
\[((j\not=j'\rightarrow t_j\not=t_{j'})\wedge
\bigwedge_{j<n} (q_j\forces t_j\in\ddot{b}))).\]
Now the claim is true by a fact about forcing (see [K, pp.201]).

\medskip

{\bf Claim 4.2}\quad 
Let $\eta\in\omega_1$ and let $q\in {\Bbb Q}$. There exists a $\nu
\leqslant\omega_1$, a maximal antichain 
$\langle p_{\alpha}:\alpha<\nu\rangle$ of $\Bbb P$, two decreasing
sequences $\langle q^j_{\alpha}:\alpha<\nu\rangle$, $j=0,1$, in $\Bbb Q$
and an increasing sequence $\langle\eta_{\alpha}:\alpha<\nu\rangle$
in $\omega_1$ such that $q^0_0, q^1_0 <q$, $\eta_0>\eta$ and for any
$\alpha<\nu$
\[p_{\alpha}\forces((\exists t_0,t_1\in\dot{T}_{\eta_{\alpha}}) 
(t_0\not=t_1\wedge (q^0_{\alpha}\forces t_0\in\ddot{b})
\wedge (q^1_{\alpha}\forces t_1\in\ddot{b}))).\]

Proof of Claim 4.2: \quad
We define those sequences inductively on $\alpha$.
First let's fix an enumeration of $\Bbb P$ in order type 
$\zeta\leqslant\omega_1$, say, ${\Bbb P}=\{x_{\gamma}:\gamma<\zeta\}$. 
For $\alpha=0$ we apply Claim 4.1 for $p=1_{\Bbb P}$
and $n=2$ to obtain $p_0, q^0_0, q^1_0$ and $\eta_0$.
Let $\alpha$ be a countable ordinal. Suppose we have found 
$\langle p_{\beta}:\beta<\alpha\rangle$,
$\langle q^0_{\beta}:\beta<\alpha\rangle$,
$\langle q^1_{\beta}:\beta<\alpha\rangle$
and $\langle\eta_{\beta}:\beta<\alpha\rangle$. If  
$\langle p_{\beta}:\beta<\alpha\rangle$ is already a maximal antichain
in $\Bbb P$, then we stop and let $\nu=\alpha$. Otherwise 
choose a smallest $\gamma<\zeta$ such that $x_{\gamma}$
is incompatible with all $p_{\beta}$'s
for $\beta<\alpha$. Pick $q^j\in {\Bbb Q}$ which are lower bounds of 
$\langle q^j_{\beta}:\beta<\alpha\rangle$ for $j=0,1$, respectively,
and pick $\eta'\in\omega_1$ which is an upper bound of
$\langle\eta_{\beta}:\beta<\alpha\rangle$. 
By applying Claim 4.1 twice we can find 
\[p'\leqslant x_{\gamma},\,\, q^0_0, q^0_1\leqslant q^0,
\,\, q^1_0, q^1_1\leqslant q^1,
\,\,\dot{t}^0_0, \dot{t}^0_1, \dot{t}^1_0, \dot{t}^1_1
\mbox{ and }\eta_{\alpha}>\eta'\] such that
\[p'\forces (\dot{t}^0_0,\dot{t}^0_1\in\dot{T}_{\eta_{\alpha}}
\wedge \dot{t}^0_0\not=\dot{t}^0_1\wedge (q^0_0\forces\dot{t}^0_0\in\ddot{b})
\wedge (q^0_1\forces\dot{t}^0_1\in\ddot{b}))\] and
\[p'\forces (\dot{t}^1_0,\dot{t}^1_1\in\dot{T}_{\eta_{\alpha}}\wedge 
\dot{t}^1_0\not=\dot{t}^1_1\wedge (q^1_0\forces\dot{t}^1_0\in\ddot{b})
\wedge (q^1_1\forces\dot{t}^1_1\in\ddot{b})).\]
If $p'\forces\dot{t}^0_0\not=\dot{t}^1_0$, then let $p_{\alpha}=p'$,
$q^0_{\alpha}=q^0_0$ and $q^1_{\alpha}=q^1_0$. 
Otherwise we can find a $p_{\alpha}<p'$
such that $p_{\alpha}\forces\dot{t}^0_0\not=\dot{t}^1_1$. Then let 
$q^0_{\alpha}=q^0_0$ and $q^1_{\alpha}=q^1_1$.
If for any countable $\alpha$, the set 
$\{p_{\beta}\in {\Bbb P}:\beta<\alpha\}$
has never been a maximal antichain, then the set 
$\{p_{\beta}\in {\Bbb P}:\beta<\omega_1\}$
must be a maximal antichain of $\Bbb P$
by the choice of $p_{\beta}$'s according to
the fixed enumeration of ${\Bbb P}=\{x_{\gamma}:\gamma<\zeta=\omega_1\}$.
In this case we choose $\nu=\omega_1$.

\medskip

The lemma follows from the construction.
Let $n\in\omega$, $\delta_n=\omega_1\cap N_n$ and let 
$\delta=\bigcup_{n\in\omega}\delta_n$.
For each $s\in 2^n$ we construct, in  $N_n$, a maximal 
antichain $\langle p^s_{\alpha}:\alpha<\nu_s\rangle$ 
of $\Bbb P$, two decreasing sequences $\langle q^{s\hat{\;}j}_{\alpha}:
\alpha<\nu_s\rangle$ for $j=0,1$,
and an increasing sequence $\langle\eta^s_{\alpha}:\alpha<\nu_s\rangle$
in $\delta_n$ such that $\nu_s\leqslant \delta_n$,
$q^{s\hat{\;}j}_0$ are lower bounds of $\langle q^s_{\alpha}: 
\alpha<\nu_{s\res n-1}\rangle$ for $j=0,1$,
$\eta^s_0=\delta^{n-1}$ and
\[p^s_{\alpha}\forces ((\exists t_0, t_1\in\dot{T}_{\eta^s_{\alpha}})
(t_0\not=t_1\wedge
(q^{s\hat{\;}0}_{\alpha}\forces t_0\in\ddot{b})\wedge 
(q^{s\hat{\;}1}_{\alpha}\forces t_1\in\ddot{b}))).\]
Each step of the construction uses Claim 4.2 relative to $N_n$ for
some $n\in\omega$. We can choose 
$q^{s\hat{\;}0}_0$ and $q^{s\hat{\;}1}_0$ to be 
lower bounds of $\langle q^s_{\alpha}: \alpha<\nu_{s\res n-1}\rangle$
because $\langle q^s_{\alpha}: \alpha<\nu_{s\res n-1}\rangle$ is constructed
in $N_{n-1}$ and hence, is countable in $N_n$. Here we use the fact
$N_{n-1}\in N_n$. 

Let $\bar{p}\leqslant 1_{\Bbb P}$ be $({\Bbb P}, N_n)$-generic for
every $n\in\omega$. 
Since $\Bbb Q$ is $\omega_1$-closed in $V$, for every $f\in 2^{\omega}$ there
is a $q_f$ which is a lower bound of $\langle q^{f\res n}_0:n\in\omega\rangle$.
Let $G\subseteq {\Bbb P}$ be a $V$-generic filter such that 
$\bar{p}\in G$. We claim that $T_{\delta}$ is uncountable
in $V[G]$. This contradicts that $T$ is an $\omega_1$-tree in $V^{\Bbb P}$.
Notice that $2^{\omega}\cap V$ is uncountable in $V[G]$. 
In $V[G]$ for each $f\in 2^{\omega}\cap V$ there is a 
$q'_f\leqslant q_f$ and a $t_f\in T_{\delta}$
such that $q'_f\forces t_f\in\dot{b}$.
Suppose $f,g\in 2^{\omega}\cap V$ are different and 
$n=\min\{i\in\omega:f(i)\not=g(i)\}$. If $t_f=t_g$, then there is a $p\in
G$, $p\leqslant\bar{p}$ such that 
\[p\forces((\exists t\in\dot{T}_{\delta}) ((q'_f\forces t\in\ddot{b})
\wedge (q'_g\forces t\in\ddot{b}))).\]
Suppose  $f\!\res\! n =s=g\!\res\! n$, $f(n)=0$ and $g(n)=1$.
Since $p$ is $({\Bbb P}, N_n)$-generic and $p\in G$, there is a 
$p^s_{\alpha}\in G$ for some $\alpha\leqslant\nu_s$.
Let $p'\leqslant p,p^s_{\alpha}$. Then 
\[p'\forces((\exists t_0,t_1\in\dot{T}_{\eta^s_{\alpha}}) 
(t_0\not=t_1\wedge (q'_f\forces t_0\in\ddot{b})\wedge
(q'_g\forces t_1\in\ddot{b}))).\]
But this contradicts the following: 
\[p'\forces(\dot{t}_0\in\dot{T}_{\eta^s_{\alpha}}\wedge\dot{t}\in
\dot{T}_{\delta}\wedge (q'_f\forces\dot{t}_0,\dot{t}\in\ddot{b})
\rightarrow \dot{t}_0\leqslant\dot{t}),\] 
\[p'\forces(\dot{t}_1\in\dot{T}_{\eta^s_{\alpha}}\wedge \dot{t}\in
\dot{T}_{\delta}\wedge (q'_f\forces\dot{t}_1,\dot{t}\in\ddot{b})
\rightarrow\dot{t}_1\leqslant\dot{t}),\] and
\[p'\forces (\dot{t}_0,\dot{t}_1\in\dot{T}_{\eta^s_{\alpha}}\wedge
\dot{t}_0\leqslant\dot{t}\wedge \dot{t}_1
\leqslant\dot{t}\rightarrow \dot{t}_0=\dot{t}_1).\]
Hence in $V[G]$ different $f$'s in $2^{\omega}\cap V$ correspond to different
$t_f$'s in $T_{\delta}$. Therefore $T_{\delta}$ is uncountable.
\quad $\Box$

\medskip

A forcing notion $\Bbb P$ is called $\omega$-proper if for any
$\omega$-sequence $\langle N_n:n\in\omega\rangle$ of countable
elementary submodels of $H(\lambda)$ such that $N_n\in N_{n+1}$
for every $n\in\omega$ and ${\Bbb P}\in N_0$, for any 
$p\in {\Bbb P}\cap N_0$ there is a $\bar{p}\leqslant p$ such that
$\bar{p}$ is $({\Bbb P},N_n)$-generic for every $n\in\omega$.
Let $S$ be a stationary subset of $\omega_1$. A forcing notion
$\Bbb P$ is called $S$-proper if for any countable elementary 
submodel $N$ of $H(\lambda)$ such that ${\Bbb P}\in N$ and
$N\cap\omega_1\in S$, and for any $p\in {\Bbb P}\cap N$ there is 
a $\bar{p}\leqslant p$ such that $\bar{p}$ is 
$({\Bbb P},N)$-generic. A forcing notion ${\Bbb P}$ is called
$(S,\omega)$-proper if for any
$\omega$-sequence $\langle N_n:n\in\omega\rangle$ of countable
elementary submodels of $H(\lambda)$ such that $N_n\in N_{n+1}$
for every $n\in\omega$, $N_n\cap\omega_1\in S$ for every $n\in\omega$,
$N\cap\omega_1\in S$, where $N=\bigcup_{n\in\omega}N_n$,
and ${\Bbb P}\in N_0$, for any 
$p\in {\Bbb P}\cap N_0$ there is a $\bar{p}\leqslant p$ such that
$\bar{p}$ is $({\Bbb P},N_n)$-generic for every $n\in\omega$.

\begin{theorem}

Let $S$ be a stationary subset of $\omega_1$ and
let $\Bbb P$ be an $(S,\omega)$-proper forcing 
notion of size $\leqslant\omega_1$ 
in $\bar{V}$. Then there are no Kurepa trees in $\bar{V}^{\Bbb P}$.

\end{theorem}

\noindent {\bf Proof:}\quad
Choose an $\eta<\kappa$ such that $S$ and $\Bbb P$ 
are in $V^{Lv(\eta,\omega_1)}$. Then \[\bar{V}^{\Bbb P}=
V^{(Lv(\eta,\omega_1)*\dot{\Bbb P})\times Lv(\kappa\smallsetminus\eta,
\omega_1)}\] and $Lv(\kappa\smallsetminus\eta,\omega_1)$ is
$\omega_1$-closed in $V^{Lv(\eta,\omega_1)}$.
By Lemma 4 it suffices to show that $\Bbb P$
satisfies property ($\dagger$) in $V^{Lv(\eta,\omega_1)}$.
Working in $V^{Lv(\eta,\omega_1)}$. Let $x\in H(\lambda)$. Since
$S$ is also stationary in $V^{Lv(\eta,\omega_1)}$,
we can choose a sequence $\langle N_n:n\in\omega\rangle$ of
countable elementary submodels of $H(\lambda)$ 
such that $N_n\in N_{n+1}$, $\{{\Bbb P}, x\}\subseteq N_0$ 
and $N_n\cap\omega_1\in S$ for every $n\in\omega$. 
Since the forcing
$Lv(\kappa\smallsetminus\eta,\omega_1)$ is countably closed, then
we can choose a decreasing sequence $\langle q_n:n\in\omega\rangle$
in $Lv(\kappa\smallsetminus\eta,\omega_1)$ such that $q_n$
is a $(Lv(\kappa\smallsetminus\eta,\omega_1), N_n)$-master condition
($q$ is a $({\Bbb Q},N)$-master condition iff for every dense
open subset $D$ of $\Bbb Q$ there exists a $d\in D$ such that
$q\leqslant d$).
Let $q$ be a lower bound of $\langle q_n:n\in\omega\rangle$.
Let $G\subseteq Lv(\kappa\smallsetminus\eta,\omega_1)$ be
$V^{Lv(\eta,\omega_1)}$-generic such that $q\in G$. By Lemma 1
every $N_n[G]$ is a countable elementary submodel of 
$(H(\lambda))^{\bar{V}}$. It is also easy to see that
$\{{\Bbb P},x\}\subseteq N_0[G]$. Now we have $N_n[G]\in N_{n+1}[G]$
and $N_n[G]\cap\omega_1\in S$ because $q\forces (N_n=N_n[\dot{G}])$.

Pick a $p\in {\Bbb P}\cap N_0$.
Since $\Bbb P$ is $(S,\omega)$-proper in $\bar{V}$, there exists
a $\bar{p}\leqslant p$ such that $\bar{p}$ is 
$({\Bbb P},N_n[G])$-generic for every $n\in\omega$. It is easy to
see that $\bar{p}$ is also $({\Bbb P},N_n)$-generic because
a maximal antichain of $\Bbb P$ in $N_n$ is also a 
maximal antichain in $N_n[G]$. This shows that $\Bbb P$
satisfies property ($\dagger$) in $V^{Lv(\eta,\omega_1)}$.
\quad $\Box$

\medskip

\noindent {\bf Remarks}\quad
(1) If $\Bbb P$ satisfies Baumgartner's axiom A, 
then $\Bbb P$ is $\omega$-proper or $(\omega_1,\omega)$-proper.
Hence forcing with a forcing notion of size 
$\leqslant\omega_1$ satisfying axiom A in $\bar{V}$
does not create Kurepa trees. Notice also that all 
c.c.c. forcing notions, $\omega_1$-closed
forcing notions and the forcing notions of tree type 
such as Sack's forcing, Laver forcing,
Miller forcing, etc. satisfy axiom A.

(2) The idea of the proof of Lemma 4 is originally from [D].
A version of Theorem 5 for axiom A forcing was proved in [J].

(3) The $\omega$-properness implies the $(S,\omega)$-properness
and the $(S,\omega)$-properness implies the property $(\dagger)$.

\bigskip

Now we prove the results about some non-$(S,\omega)$-proper forcing
notions.

The existence of a Kurepa tree implies that there are no countably
complete, $\aleph_2$-saturated ideals on $\omega_1$. 
Therefore, one can destroy all those ideals
by creating a generic Kurepa tree [V]. But one don't have to
create Kurepa trees for this purpose. Baumgartner and Taylor [BT]
proved that adding a club subset of $\omega_1$ by finite
conditions destroys all countably complete, $\aleph_2$-saturated 
ideals on $\omega_1$.
The forcing notion for adding a club subset of $\omega_1$ by finite
conditions has size $\leqslant\omega_1$ and is proper but not
$(S,\omega)$-proper for any stationary subset $S$ of $\omega_1$..
We are going to prove next that this forcing notion and some other
similar forcing notions do not create
Kurepa trees if our ground model is the L\'{e}vy model $\bar{V}$.
Notice also that the ideal of nonstationary subsets of $\omega_1$
could be $\aleph_2$-saturated in the L\'{e}vy model obtained by
collapsing a supercompact cardinal down to $\omega_2$ [FMS].
As a corollary we can have a ground model $\bar{V}$ which contains 
countably complete, $\aleph_2$-saturated ideals on $\omega_1$
such that forcing with some small proper forcing notion $\Bbb P$
in $\bar{V}$
destroys all countably complete, $\omega_2$-saturated ideals on 
$\omega_1$ without creating Kurepa trees.

We first define a property of forcing notions which is satisfied by
the forcing notion for adding a club subset of $\omega_1$
by finite conditions.

\begin{definition}

A forcing notion $\Bbb P$ is said to satisfy property (\#) if
for any $x\in H(\lambda)$
there exists a countable elementary submodel $N$ of $H(\lambda)$ such that
$\{{\Bbb P},x\}\subseteq N$ and for any 
$p_0\in {\Bbb P}\cap N$ there exists a 
$\bar{p}\leqslant p_0$, $\bar{p}$ is $({\Bbb P}, N)$-generic,
and there exists a countable subset $C$ of $\Bbb P$ such that
for any $\bar{p}'\leqslant\bar{p}$ there is a $c\in C$ and a 
$p'\in {\Bbb P}\cap N$, $p'\leqslant p_0$ such that

(1) for any dense open subset $D$ of ${\Bbb P}$ below $p'$ in $N$
there is an $d\in D\cap N$ such that $d$ is compatible with $c$, and

(2) for any $r\in {\Bbb P}\cap N$ and $r\leqslant p'$,
$r$ is compatible with $c$ implies $r$ is compatible with $\bar{p}'$.

Let's call the pair $(p',c)$ a related pair corresponding to
$\bar{p}'$.

\end{definition}

\begin{example}

Following three examples are the forcing notions which satisfy
property (\#).

\end{example}

(1) Let \[
\begin{array}{ll}
{\Bbb P}= & \{p\subseteq\omega_1\times\omega_1:
p\mbox{ is a finite function which can be extended to}\\
 &\mbox{an increasing continuous function from }\omega_1\mbox{ to }\omega_1.\}
\end{array} \] and
let $\Bbb P$ be ordered by reverse inclusion. 
$\Bbb P$ is one of the simplest proper forcing notion which
does not satisfy axiom A [B2]. 
Forcing with $\Bbb P$ creates a generic
club subset of $\omega_1$ and destroys all $\aleph_2$-saturated
ideals on $\omega_1$ [BT]. It is easy to see that
$\Bbb P$ satisfies property (\#) defined above. For any 
$x\in H(\lambda)$ we can choose a countable elementary submodel $N$
of $H(\lambda)$ such that $\{{\Bbb P}, x\}\subseteq N$ and
$N\cap\omega_1=\delta$ is an indecomposable ordinal. For any $p_0\in 
{\Bbb P}\cap N$ let $\bar{p}=p_0\cup (\delta,\delta)$ and let
$C=\{\bar{p}\}$. Then for any $\bar{p}'\leqslant\bar{p}$ there is a 
$p'=\bar{p}'\!\res\!\delta$ and a $c=\bar{p}\in C$ such that all requirements
for the definition of property (\#) are satisfied.

\bigskip

(2) Let $S$ be a stationary subset of $\omega_1$.
If we define 
\[
\begin{array}{ll}
{\Bbb P}_S= & \{p: p\mbox{ is a finite function such that there is an
increasing continuous}\\
 & \mbox{ function }f\mbox{ from some countable ordinal
to }S\mbox{ such that }p\subseteq f.\}
\end{array}
\]
 and let ${\Bbb P}_S$ be
ordered by reverse inclusion, then ${\Bbb P}_S$ is $S$-proper [B2].
Forcing with ${\Bbb P}_S$ adds a club set inside $S$.
It is also easy to check that ${\Bbb P}_S$ satisfies (\#).
For any $x\in H(\lambda)$. Let $N$ be a countable elementary submodel
of $H(\lambda)$ such that $\{x,{\Bbb P}_S\}\subseteq N$,
$N\cap\omega_1=\delta$ is an indescomposable ordinal and 
$\delta\in S$. Then for any $p_0\in {\Bbb P}_S\cap N$ the element
$\bar{p}=p_0\cup\{(\delta,\delta)\}$ is $({\Bbb P}_S,N)$-generic.
Now $N$, $\bar{p}$ and $C=\{\bar{p}\}$ witness that ${\Bbb P}_S$
satisfies property (\#).

\bigskip

(3) Let $T$ and $U$ be two normal Aronszajn trees such that every
node of $T$ or $U$ has infinitely many immediate successors.
Let $\Bbb P$ be the forcing notion such that 
$p=(A_p,f_p)\in {\Bbb P}$ iff 

(a) $A_p$ is a finite subset of $\omega_1$, 

(b) $f_p$ is a finite partial isomorphism from
$T\!\res\! A_p$ into $U\!\res\! A_p$, 

(c) $dom(f_p)$ is a subtree of
$T\!\res\! A_p$ in which every branch has cardinality $|A_p|$.

\noindent $\Bbb P$ is ordered by $p\leqslant q$ iff $A_p\supseteq A_q$
and $f_p\supseteq f_q$. $\Bbb P$ is proper [T]. $\Bbb P$ is used
in [AS] for generating a club isomorphism from
$T$ to $U$. For any $x\in H(\lambda)$, for any countable elementary
submodel $N$ of $H(\lambda)$ such that $\{{\Bbb P},x\}\subseteq N$
and for any $p_0\in {\Bbb P}\cap N$, let $\delta=N\cap\omega_1$,
let $A_{\bar{p}}=A_{p_0}\cup\{\delta\}$ and 
let $f_{\bar{p}}$ be any extension of $f_{p_0}$ such that 
$T_{\delta}\cap dom(f_{\bar{p}})\not=\emptyset$. Then
$\bar{p}=(A_{\bar{p}},f_{\bar{p}})$ is a $({\Bbb P},N)$-generic
condition. Let \[C=\{d:d\mbox{ is a finite isomorphism from }T_{\delta}
\mbox{ to }U_{\delta}\}.\] Then $C$ is countable. For any $\bar{p}'\leqslant \bar{p}$
let $c=(f_{\bar{p}'}\!\res\!\{\delta\}) \in C$, let $\alpha<\delta$, $\alpha>
\max (A_{\bar{p}'}\cap\delta)$ and 
\[g_{\alpha}=\{(t,u)\in T_{\alpha}\times
U_{\alpha}:(\exists (t',u')\in (f_{\bar{p}'}\!\res\!\{\delta\})) 
(t<t'\wedge u<u')\}\] be such that $g_{\alpha}$ and 
$f_{\bar{p}'}\!\res\!\{\delta\}$
have same cardinality, let $A_{p'}=(A_{\bar{p}'}\cap\delta)\cup\{\alpha\}$, 
let $f_{p'}=(f_{\bar{p}'}\!\res\! (A_{\bar{p}'}\cap\delta))\cup g_{\alpha}$, 
and let $p'=(A_{p'},f_{p'})$. 
Then $(p',c)$ is a related pair corresponding to $\bar{p}'$ [AS] and
$N,\bar{p},C$ witness that $\Bbb P$ satisfies 
property (\#).  For any stationary set $S$ we can also define
an $S$-proper version of this forcing notion.

\begin{lemma}

Let $V$ be a model. Let $\Bbb P$ and $\Bbb Q$ be two forcing notions
in $V$ such that $\Bbb P$ has size $\leqslant\omega_1$ and satisfies property 
(\#), and $\Bbb Q$ is $\omega_1$-closed (in $V$). 
Suppose $T$ is an $\omega_1$-tree in $V^{\Bbb P}$.
Then $T$ has no branches which are in $V^{{\Bbb P}\times {\Bbb Q}}$ but
not in $V^{\Bbb P}$.

\end{lemma}

\noindent {\bf Proof:}\quad
Suppose, towards a contradiction, that there is a branch $b$ of $T$ in 
$V^{{\Bbb P}\times {\Bbb Q}}\smallsetminus V^{\Bbb P}$.
Without loss of generality, we assume that
\[\forces_{\Bbb P}\,\forces_{\Bbb Q} 
(\ddot{b}\mbox{ is a branch of }\dot{T}\mbox{ in }
V^{{\Bbb P}\times {\Bbb Q}}\smallsetminus V^{\Bbb P}).\]
Following the definition of property (\#), we can find
a countable elementary submodel $N$ of $H(\lambda)$ 
such that $\{{\Bbb P}, {\Bbb Q},\dot{T},\ddot{b}\}\subseteq N$,
a $\bar{p}\leqslant 1_{\Bbb P}$ which is 
$({\Bbb P},N)$-generic and a countable set
$C\subseteq {\Bbb P}$ such that
$N$, $\bar{p}$ and $C$ witness that $\Bbb P$ satisfies
property (\#). Let $\langle (p_i,c_i):i\in\omega\rangle$
be a listing of all related pairs in $({\Bbb P}\cap N)\times C$
with infinite repetition, {\em i.e.} every related pair $(p,c)$
in $({\Bbb P}\cap N)\times C$ occurs infinitely ofter in
the sequence.

We construct now, in $V$,
a set $\{q_s\in {\Bbb Q}\cap N :s\in 2^{<\omega}\}$ and
an increasing sequence $\langle\delta_n:n\in\omega\rangle$ such that

(1) $s\subseteq t$ implies $q_t\leqslant q_s$,

(2) $\delta_n\in\delta=N\cap\omega_1$,

(3) for every $n\in\omega$ there is a $p'\in {\Bbb P}\cap N, p'\leqslant
p_n$ such that $p'$ is compatible with $c_n$, and 
\[p'\forces ((\exists\{t_s:s\in 2^n\}\subseteq\dot{T}_{\delta_n})
((s\not=s'\rightarrow t_s\not=t_{s'})\wedge \bigwedge_{s\in 2^n}
(q_s\forces t_s\in\ddot{b}))).\]

The lemma follows from the construction. Let $G\subseteq {\Bbb P}$ be
a $V$-generic filter and $\bar{p}\in G$. We want to show that
\[V[G]\models T_{\delta}\mbox{ is uncountable.}\] For any $f\in 2^{\omega}\cap V$
let $q_f\in {\Bbb Q}$ be a lower bound of the set 
$\{q_{f\res n}:n\in\omega\}$ such that there is
a $t_f\in T_{\delta}$ such that
$q_f\forces t_f\in\dot{b}$. Suppose $T_{\delta}$ is countable. Then
there are $f,g\in 2^{\omega}\cap V$ such that $t_f=t_g$. Let
$\dot{t}_f$, $\dot{t}_g$ be $\Bbb P$-names for
$t_f,t_g$ and let $\bar{p}'\leqslant\bar{p}$ be such that 
\[\bar{p}'\forces (\dot{t}_f=\dot{t}_g\wedge (q_f\forces \dot{t}_f\in\ddot{b})
\wedge (q_g\forces\dot{t}_g\in\ddot{b})).\]
Let $m=\min\{i\in\omega:f(i)\not=g(i)\}$.
By the definition of property (\#) we can find a related pair $(p,c)$ 
corresponding to $\bar{p}'$. Choose an $n\in\omega$
such that $n\geqslant m$ and $(p,c)=(p_n,c_n)$. Since (1) of 
Definition 6 is true, there is a $p'\in {\Bbb P}\cap N$
such that $p'\leqslant p$, $p'$ is compatible with $c_n$ and
\[p'\forces ((\exists\{t_s:s\in 2^n\}\subseteq\dot{T}_{\delta_n})
((s\not=s'\rightarrow t_s\not=t_{s'})\wedge
\bigwedge_{s\in 2^n}(q_s\forces t_s\in\ddot{b}))).\]
Since $q_f\leqslant q_{f\res n}$ and $q_g\leqslant q_{g\res n}$, then 
\[\bar{p}'\forces ((\exists t_0,t_1\in\dot{T}_{\delta_n}) 
(t_0\not=t_1\wedge (q_f\forces t_0\in\ddot{b})\wedge (q_g\forces t_1\in
\ddot{b}))).\]
But also \[\bar{p}'\forces((\exists t\in\dot{T}_{\delta}) ((q_f\forces
t\in\ddot{b})\wedge (q_g\forces t\in\ddot{b}))).\]
By the fact that any two nodes in $T_{\delta_n}$ which are
below a node in $T_{\delta}$ must be same, and that $p'$ is compatible
with $\bar{p}'$, we have a contradiction.

Now let's inductively construct $\{\delta_i:i\in\omega\}$ and
$\{q_s:s\in 2^{<\omega}\}$. 
Suppose we have had $\{q_s:s\in 2^{\leqslant n}\}$
and $\{\delta_i:i\leqslant n\}$. let $D\subseteq {\Bbb P}$ be such that
$r\in D$ iff 

(1) $r\leqslant p_n$ (recall that $(p_n,c_n)$ is in the enumeration of
all related pairs in $({\Bbb P}\cap N)\times C$),

(2) there exists $\eta >\delta_n$ and there exists 
$\{q_s\leqslant q_{s\res n}:s\in 2^{n+1}\}$ such that
\[ r\forces ((\exists \{t_s:s\in 2^{n+1}\}\subseteq\dot{T}_{\eta})
((s\not=s'\rightarrow t_s\not=t_{s'})\wedge
\bigwedge_{s\in 2^{n+1}}(q_s\forces t_s\in\ddot{b}))).\]
It is easy to see that $D$ is open and $D\in N$.

\medskip

{\bf Claim 8.1}\quad $D$ is dense below $p_n$.

Proof of Claim 8.1:\quad Suppose $r_0\leqslant p_n$. It suffices to
show that there is an $r\leqslant r_0$ such that $r\in D$.
Applying Claim 4.1, for any $s\in 2^n$ we can 
find $r_s\leqslant r_0$, $\eta_s>\delta_n$
and $\{q^s_j\leqslant q_s:j<2^{n+1}\}$ such that
\[ r_s\forces ((\exists \{t_j:j<2^{n+1}\}\subseteq\dot{T}_{\eta_s})
((j\not=j'\rightarrow t_j\not=t_{j'})\wedge \bigwedge_{j<2^{n+1}}
(q^s_j\forces t_j\in\ddot{b}))).\]
Let $\{s_i:i<2^n\}$ be an enumeration of $2^n$. By applying
Claim 4.1 $2^n$ times as above we obtained $r_0\geqslant r_{s_0}
\geqslant r_{s_1}\geqslant\ldots r_{s_{2^n-1}}$ such that above
arguments are true for any $s\in 2^n$. Pick 
$\eta=\max\{\eta_s:s\in 2^n\}$. Then we extend $r_{s_{2^n-1}}$
to $r'$, and extend $q^s_j$ to $\bar{q}^s_j$ for every such $s$ and $j$
such that for each $s\in 2^n$
\[ r'\forces ((\exists \{t_j:j<2^{n+1}\}\subseteq\dot{T}_{\eta})
((j\not=j'\rightarrow t_j\not=t_{j'})\wedge\bigwedge_{j<2^{n+1}}
(\bar{q}^s_j\forces t_j\in\ddot{b}))).\]
Now applying an argument in Claim 4.2 repeatedly 
we can choose $\{q_{s\hat{\;}0},q_{s\hat{\;}1}\}\subseteq$\break\hfill
$\{\bar{q}^s_j:j<2^{n+1}\}$ for every $s\in 2^n$ and extend $r'$ to $r''$
such that 
\[ r''\forces ((\exists \{t_s:s\in 2^{n+1}\}\subseteq\dot{T}_{\eta})
((s\not=s'\rightarrow t_s\not=t_{s'})\wedge
\bigwedge_{s\in 2^{n+1}}(q_s\forces t_s\in\ddot{b}))).\]
This showed that $D$ is dense below $p_n$.

\medskip

Notice that since $N$ is elementary, then $\eta$ exists in $N$ and
all those $q_s$' for $s\in 2^{n+1}$ exist in $N$.
Choose $r\in D$ such that $r,c_n$ are compatible and let
$\delta_{n+1}$ be correspondent $\eta$. 
This ends the construction. \quad $\Box$

\begin{theorem}

If ${\Bbb P}$ in $\bar{V}$ is a forcing notion defined in 
(1), (2) or (3) of Examples 7,
then forcing with $\Bbb P$ does not create any Kurepa trees.

\end{theorem}

\noindent {\bf Proof:}\quad
Suppose $T$ is a Kurepa tree in $\bar{V}^{\Bbb P}$.
Let $\eta<\kappa$ be such that ${\Bbb P}, T\in V^{Lv(\eta,\omega_1)}$.
Since the definition of $\Bbb P$ is absolute between
$\bar{V}$ and $V^{Lv(\eta,\omega_1)}$, then $\Bbb P$ satisfies
property (\#) in $V^{Lv(\eta,\omega_1)}$. Since
$T$ has less than $\kappa$ branches in 
$V^{Lv(\eta,\omega)*\dot{\Bbb P}}$, there exist branches of $T$ in
$\bar{V}^{\Bbb P}$ which are not in $V^{Lv(\eta,\omega_1)*\dot{\Bbb P}}$.
This contradicts Lemma 8. \quad $\Box$

\bigskip

\noindent {\bf Remark:}\quad The forcing notions in Examples 7, (1), (2)
and (3) are not $(S,\omega)$-proper for any stationary $S$.

\section{Creating Kurepa Trees By a Small Forcing Is Easy}

In this section we construct a model of \ch\, plus no Kurepa
trees, in which there is an $\omega$-distributive Aronszajn tree $T$ 
such that forcing with $T$ does create a Kurepa tree in the generic extension.

Let $V$ be a model and $\kappa$ be a strongly inaccessible 
cardinal in $V$. Let $\cal T$ be the set of all countable normal
trees. Given a set $A$ and a cardinal $\lambda$. Let 
$[A]^{<\lambda} = \{S\subseteq A:|S|<\lambda\}$ and
$[A]^{\leqslant\lambda} = \{S\subseteq A:|S|\leqslant\lambda\}$.
We define a forcing notion $\Bbb P$ as following:

\begin{definition}

$p$ is a condition in $\Bbb P$ iff 
\[p=\langle\alpha_p,t_p,k_p,U_p,B_p,F_p\rangle\]
where

(a) $\alpha_p\in\omega_1$,

(b) $t_p\in {\cal T}$ and $ht(t_p)=\alpha_p+1$,

(c) $k_p$ is a function from $t_p$ to $\cal T$ such that for any 
$x\in t_p$, $ht(k_p(x))=ht(x)+1$,
and for any $x,y\in t_p$, $x<y$ implies $k_p(x)\leqslant_{end}k_p(y)$,

(d) $U_p\in [\kappa]^{\leqslant\omega_1}$,

(e) $B_p=\{b^p_{\gamma}:\gamma\in U_p\}$ where $b^p_{\gamma}$ is a
function from $t_p\!\res\! (\beta^p_{\gamma}+1)$ to $\omega_1^{<\omega_1}$
for some $\beta^p_{\gamma}\leqslant\alpha_p$ such that for any 
$x\in t_p\!\res\! (\beta^p_{\gamma}+1)$,
$b^p_{\gamma}(x)\in (k_p(x))_{ht(x)}$ and for any $x,y\in t_p\!\res\!
(\beta^p_{\gamma})$, $x\leqslant y$ implies $b^p_{\gamma}(x)\leqslant
b^p_{\gamma}(y)$,

(f) $F_p=\{f^p_{\gamma}:\gamma\in U_p\}$ where $f^p_{\gamma}$ is a
function from $\delta^p_{\gamma}$ to $\gamma$ for some
$\delta^p_{\gamma}\leqslant\alpha_p$,

(g) for any $x\in t_p\!\res\!\alpha_p$, 
for any finite $U_0\subseteq U_p$ and for any
$\epsilon$ such that $ht(x)<\epsilon\leqslant\alpha_p$, there exists
an $x'\in (t_p)_{\epsilon}$ such that $x'>x$ and for any
$\gamma_1,\gamma_2\in U_0$ either one of $\beta^p_{\gamma_1},
\beta^p_{\gamma_2}$ is less than $\epsilon$ or
$b^p_{\gamma_1}(x)=b^p_{\gamma_2}(x)$ implies
$b^p_{\gamma_1}(x')=b^p_{\gamma_2}(x')$.

\end{definition}

In the condition (g) of the definition we call $x'$ a conservative
extension of $x$ at level $\epsilon$ with respect to $U_0$ (or
with respect to $\{b^p_{\gamma}:\gamma\in U_0\}$).

Generally we have the following notation.
Suppose $t\in {\cal T}$ and $B$ is a set of functions such that
for each $b\in B$ there is a $\beta_b\leqslant ht(t)$ such that
$domain(b)=t\!\res\!\beta$. We say $t$ is consistent with respect to $B$ if 
for any $x\in t\!\res\! ht(t)$, for any finite $B_0\subseteq B$ and for any
$\epsilon$ such that $ht(x)<\epsilon\leqslant ht(t)$, there exists
an $x'\in t_{\epsilon}$ such that $x'>x$ and for any
$b_1,b_2\in B_0$ either one of $\beta_{b_1},
\beta_{b_2}$ is less than $\epsilon$ or
$b_1(x)=b_2(x)$ implies $b_1(x')=b_2(x')$. So $p\in {\Bbb P}$ implies
that $t_p$ is consistent with respect to $B_p$.

\medskip

For any $p,q\in {\Bbb P}$ we define the order of $\Bbb P$ by letting
$p\leqslant q$ iff

(1) $\alpha_q\leqslant\alpha_p$, $t_q\leqslant_{end}t_p$,
$k_q\subseteq k_p$ and $U_q\subseteq U_p$, 

(2) for any $\gamma\in U_q$, $b^q_{\gamma}\subseteq b^p_{\gamma}$ and
$f^q_{\gamma}\subseteq f^p_{\gamma}$,

(3) $\{\gamma\in U_q:\beta^p_{\gamma}>\beta^q_{\gamma}\}$ is at most
countable, 

(4) $\{\gamma\in U_q:\delta^p_{\gamma}>\delta^q_{\gamma}\}$ is at most
countable.

\bigskip

\noindent {\bf Remarks:}\quad In the definition of $\Bbb P$
the part $t_p$ is used for creating an $\omega$-distributive Aronszajn
tree $T$. The part $k_p$ is used for creating a $T$-name of an $\omega_1$-tree
$K$. The part $B_p$ is used for adding $\kappa$ branches to $K$ so that
$K$ becomes a Kurepa tree in the generic extension by forcing with $T$.
The part $F_p$ is used for collapsing all cardinals between
$\omega_1$ and $\kappa$.

\bigskip

For any $\epsilon\in\omega_1$, $\gamma\in\kappa$ and $\eta\in\gamma$,
let
\[D^1_{\epsilon} = \{p\in {\Bbb P}:\alpha_p\geqslant\epsilon\},\]
\[D^2_{\gamma} = \{p\in {\Bbb P}:\gamma\in U_p\},\]
\[D^3_{\eta,\gamma} = \{p\in {\Bbb P}:\gamma\in U_p\mbox{ and }
\eta\in range(f^p_{\gamma})\},\]
\[D^4_{\epsilon,\gamma} = \{p\in {\Bbb P}:\gamma\in U_p\mbox{ and }
\beta^p_{\gamma}\geqslant\epsilon\}.\]

\begin{lemma}

The sets $D^1_{\epsilon}$, $D^2_{\gamma}$, $D^3_{\eta,\gamma}$ and
$D^4_{\epsilon,\gamma}$ are open dense in $\Bbb P$.

\end{lemma}

\noindent {\bf Proof:}\quad
It is easy to see that all four sets are open. Let's show they are
dense. The proofs of the denseness of the first three sets are easy.

Given $p_0\in {\Bbb P}$. We need to find a $p\leqslant p_0$ such that
$p\in D^1_{\epsilon}$. Pick an $\alpha_p\geqslant\epsilon$ and
$\alpha_p\geqslant\alpha_{p_0}$. Let $t_p\in {\cal T}$ be such that
$ht(t_p) = \alpha_p+1$ and $t_{p_0}\leqslant_{end}t_p$. Let 
$k_p:t_p\mapsto {\cal T}$ be any suitable extension of $k_{p_0}$.
Let $U_p=U_{p_0}$. For any $\gamma\in U_p$ let
$b^p_{\gamma}=b^{p_0}_{\gamma}$ and 
$f^p_{\gamma}=f^{p_0}_{\gamma}$. Then $p\leqslant p_0$ and 
$p\in D^1_{\epsilon}$.

Given $p_0\in {\Bbb P}$. We need to find a $p\leqslant p_0$ such that
$p\in D^2_{\gamma}$. If $\gamma\in U_{p_0}$, let $p=p_0$. Otherwise,
let 
\[p=\langle\alpha_{p_0},t_{p_0},k_{p_0},U_{p_0}\cup\{\gamma\},
B_{p_0}\cup\{b^p_{\gamma}\},F_{p_0}\cup\{f^p_{\gamma}\}\rangle,\] where
$b^p_{\gamma}$ and $f^p_{\gamma}$ are empty functions. Then $p\leqslant
p_0$ and $p\in D^2_{\gamma}$.

Given $p_0\in {\Bbb P}$. We need to find a $p\leqslant p_0$ such that
$p\in D^3_{\eta,\gamma}$. First, pick $p'\in D^1_{\alpha_0+1}$ such that
$p'\leqslant p_0$ and $f^{p'}_{\gamma} = f^{p_0}_{\gamma}$. Then extend
$f^{p'}_{\gamma}$ to $f^p_{\gamma}$ on $\alpha_0+1$ arbitrary except
assigning $f^p_{\gamma}(\alpha_0)
= \eta$. Let everything else keep unchanged. Then $p\leqslant p'$ and
$p\in D^3_{\eta,\gamma}$.

Proving the denseness of $D^4_{\epsilon,\gamma}$ is not trivial due to
the condition (g) of Definition 10. Given $p_0\in {\Bbb P}$. Without
loss of generality we assume that $p_0\in D^1_{\epsilon}\cap
D^2_{\gamma}$ and $\epsilon > \beta^{p_0}_{\delta}$ for all $\delta\in
U_{p_0}$. We need to find a $p\leqslant p_0$ such that $p\in
D^4_{\epsilon,\gamma}$. Choose $\alpha_p=\alpha_{p_0},t_p=t_{p_0},
k_p=k_{p_0},U_p=U_{p_0},b^p_{\delta}=b^{p_0}_{\delta}$ for all
$\delta\in U_{p_0}\smallsetminus\{\gamma\}$ and $f^p_{\delta}=
f^{p_0}_{\delta}$ for all $\delta\in U_{p_0}$. Let $\beta^p_{\gamma}
=\epsilon$. We need to extend $b^{p_0}_{\gamma}$ to $b^p_{\gamma}$ on
$t_p\!\res\! (\epsilon+1)$ such that $p\in {\Bbb P}$.

For each $x\in t_p\!\res\! (\epsilon+1)\smallsetminus t_p\!\res\!
\beta^{p_0}_{\gamma}$ and for each $\mu\leqslant\epsilon$ Let
$C_{x,\mu}$ be the cone above $x$ up to level $\mu$, {\em i.e.} 
\[C_{x,\mu}=\{y\in t_p:x<y\mbox{ and }ht(y)\leqslant\mu\}.\] 
We construct $t_0\subseteq t_1\subseteq\ldots$ with
$t_0=t_p\!\res\!\beta^{p_0}_{\gamma}$ and define $b^p_{\gamma}$ on 
$t_n$ inductively. Suppose we have had $t_n$ and $b^p_{\gamma}\!\res\! t_n$.
For any maximal node $x$ of $t_n$ we define a subset $t^n_x$ above $x$.
It will be self-clear from the construction that for any $n\in\omega$
and for any $x\in t_n$ there is a maximal node $x'$ of $t_n$ such
that $x'\geqslant x$.
Our $t_{n+1}$ will be the union of $t_n$ and those $t^n_x$'s.
Let $x$ be a maximal node of $t_n$. Let \[U_x = 
\{\beta^p_{\delta}: \delta\in U_p\smallsetminus\{\gamma\},
\,\beta^p_{\delta}>ht(x) \mbox{ and }
b^p_{\delta}(x)=b^p_{\gamma}(x) \}.\]

Case 1:\quad $U_x=\emptyset$. Let $t^n_x=\emptyset$. This means any
choice of $b^p_n$ above $x$ will not violate the condition (g).

Case 2:\quad $U_x$ has a largest element, say $\beta^p_{\delta'}$.
Let $t^n_x=C_{x,\beta^p_{\delta'}}$ and let $b^p_{\gamma}\!\res\! t^n_x
=b^p_{\delta'}\!\res\! t^n_x$.

Case 3:\quad $\bigcup U_x$ is a limit ordinal. Fix a strictly increasing
sequence $\langle\nu_{x,m}:m\in\omega\rangle$ 
of ordinals such that $\bigcup_{m\in\omega}\nu_{x,m}=\bigcup U_x$.
Let $x_0\leqslant x_1\leqslant\ldots\leqslant x_n=x$ be such that
$x_i$ is a maximal node of $t_i$ for $i=0,1,\ldots,n$. Notice that 
if $i<n$, then $\bigcup U_{x_i}\geqslant\bigcup U_x$, 
and if $\,\bigcup U_{x_i}$
is a limit ordinal, then $\langle\nu_{x_i,m}:m\in\omega\rangle$
has already been defined. Let 
\[l=\min\{i:\bigcup_{m\in\omega}\nu_{x_i,m}=
\bigcup_{m\in\omega}\nu_{x_n,m}\}\] and let
\[\bar{\nu}=\max\{\nu_{x_i,n}:l\leqslant i\leqslant n\}.\] Choose
$\delta\in U_x$ such that $\beta^p_{\delta}\geqslant\bar{\nu}$ and
let $b^p_{\gamma}\!\res\! C_{x,\beta^p_{\delta}}=b^p_{\delta}\!\res\! 
C_{x,\beta^p_{\delta}}$. Let $t^n_x = C_{x,\beta^p_{\delta}}$.
Now we take 
\[t_{n+1} = t_n\cup (\bigcup\{t^n_x:x\mbox{ is a maximal node
of }t_n.\}\] and define $b^p_{\gamma}\!\res\! t_{n+1}$ accordingly.
Let $t=\bigcup_{n\in\omega}t_n$. Notice that $t$ may not be equal to $t_p\!\res\! 
(\epsilon+1)$. But it is no problem because any extension of
$b^p_{\gamma}\!\res\! t$ to $t_p\!\res\! (\epsilon+1)$ following 
the condition (e) will not violate the condition (g). Let $b^p_{\gamma}$
be such an extension of $b^p_{\gamma}\!\res\! t$.

\medskip

{\bf Claim 11.1}\quad $p\in {\Bbb P}$.

Proof of Claim 11.1:\quad We need only to check that the condition (g)
of Definition 10 is satisfied. Pick $x\in t_p\!\res\!\epsilon$ and pick
a finite subset $U_0$ of $U_p$. Pick also an $\epsilon'$ such that
$ht(x)<\epsilon'\leqslant\epsilon$. First, we assume that 
$x\in t_n\smallsetminus t_{n-1}$ for some $n\in\omega$ (let
$t_{-1}=\emptyset$). Without loss of generality we assume that
$x$ is a maximal node of $t_n$.

Case 1:\quad Every $\beta\in U_x$ is less than $\epsilon'$. Then
the condition (g) is trivially satisfied because any conservative
extension of $x$ at level $\epsilon'$ with respect to $U_0\smallsetminus
\{\gamma\}$ is a conservative extension of $x$ with respect to
$U_0$.

Case 2:\quad There is a largest ordinal $\beta^p_{\delta'}\geqslant
\epsilon'$ in $U_x$ such that 
\[b^p_{\gamma}\!\res\! C_{x,\beta^p_{\delta'}}
=b^p_{\delta'}\!\res\! C_{x,\beta^p_{\delta'}}.\] Then a conservative extension
of $x$ at level $\epsilon'$ with respect to $(U_0\smallsetminus\{\gamma\})
\cup\{\delta'\}$ is a conservative extension of $x$ with respect to
$U_0$.

Case 3:\quad $\bigcup U_x$ is a limit ordinal greater than $\epsilon'$.
First, choose $\beta^p_{\delta'}>\epsilon'$ in $U_x$. Suppose $\nu_{x,m}
\leqslant\beta^p_{\delta'}<\nu_{x,m+1}$. Then choose a maximal node $x_1$
of $t_{n+1}$ such that $x_1$ is a conservative extension of $x$ with
respect to $U_0\cup\{\gamma,\delta'\}$. Now we have \[\bigcup
U_{x_1}\geqslant\beta^p_{\delta'}>\epsilon'.\] Notice that 
$ht(x_1)\geqslant \nu_{x,n}$. We are done if $U_{x_1}$ has a largest
ordinal. Otherwise we repeat the same procedure to get $x_2$. Eventually,
we can find an $x_k$ such that $x_k$ is a conservative extension of $x$
with respect to $U_0\cup\{\delta'\}$ and $ht(x_k)\geqslant
\nu_{x,m+1}>\epsilon'$. Let $x''\leqslant x'$ and $ht(x'')=\epsilon'$.
It is easy to see that $x''$ is a conservative extension of $x$
at level $\epsilon'$ with respect to $U_0$.

Suppose $x\not\in t$. Then $U_x=\emptyset$. So every $x'\geqslant x$,
$x'\in t_{\epsilon'}$ is a conservative extension of $x$ with
respect to $U_0$.

This ends the proof of the claim. It is easy to see that $p\in
D^4_{\epsilon,\gamma}$.\quad $\Box$

\bigskip

Next we want to prove that ${\Bbb P}$ is $\omega_1$-strategically
closed. Let $\Bbb Q$ be a forcing notion. Two players, $I$ and $I\!I$,
play a game $G({\Bbb Q})$ by $I$ choosing $p_n\in {\Bbb Q}$ and $I\!I$
choosing $q_n\in {\Bbb Q}$ alternatively 
such that \[p_0\geqslant q_0\geqslant
p_1\geqslant q_1\geqslant\ldots.\] $I\!I$ wins the game $G({\Bbb Q})$ if
and only if the sequence $\langle p_0,q_0,p_1,q_1,\ldots\rangle$ has
a lower bound in $\Bbb Q$. A forcing notion $\Bbb Q$ is called
$\omega_1$-strategically closed if $I\!I$ wins the game $G({\Bbb Q})$.
Note that any $\omega_1$-strategically closed forcing notion does
not add new countable sequences of ordinals to the generic extension.

\begin{lemma}

$\Bbb P$ is $\omega_1$-strategically closed.

\end{lemma}

\noindent {\bf Proof:}\quad
We choose $q_n$ inductively for Player $I\!I$ after Player $I$ choose
any $p_n\leqslant q_{n-1}$. Suppose $p_i,q_i$ have been chosen for
$i<n$. Let $p_n\leqslant q_{n-1}$ be any element chosen by Player $I$.
Player $I\!I$ want to choose $q_n\leqslant p_n$. Let 
\[U_n=\{\gamma\in
U_{p_n}:(\exists i<n) (\beta^{p_i}_{\gamma}\not=\beta^{q_i}_{\gamma})
\mbox{ or }(\exists i\leqslant n) (\beta^{p_i}_{\gamma}\not=
\beta^{q_{i-1}}_{\gamma})\}.\] Choose $q_n\leqslant p_n$ such that
$\alpha_{q_n}>\alpha_{p_n}$ and for any $\gamma\in U_n$,
$\beta^{q_n}_{\gamma} = \alpha_{q_n}$. This can be done by
repeating the steps countably many times used in the proof
of the denseness of $D^4_{\epsilon,\gamma}$ in Lemma 11.
This finishes the inductive step of the construction.
Let \[\alpha_q =\bigcup_{n\in\omega}\alpha_{q_n},\; t' =
\bigcup_{n\in\omega},\; k'=\bigcup_{n\in\omega}k_{q_{n}},\;
U_q=\bigcup_{n\in\omega}U_{q_n}\] and for each $\gamma\in U_q$
\[b'_{\gamma}=\bigcup\{b^{q_n}_{\gamma}:n\in\omega,\gamma\in
U_{q_n}\}\] and \[f^q_{\gamma}=\bigcup\{f^{q_n}_{\gamma}:n\in\omega,\gamma\in
U_{q_n}\}.\] We need now to add one more level on the top of $t'$
and extend $k'$ and $b'_{\gamma}$'s accordingly.
The main difficulty here is to make the condition (g) of Definition 10 true.
Remember \[U_{\omega} =\bigcup_{n\in\omega}U_n\subseteq U_p\] is the set
of all $\gamma$'s such that $\beta^{q_n}_{\gamma}$ grows for some $n$. 
The set $U_{\omega}$ is at most countable 
due to the definition of the order of
$\Bbb P$. Note that $\alpha_{q_n}$ is strictly increasing. Note also
that for each $\gamma\in U_q\smallsetminus U_{\omega}$ the sequence
\[\{b^{q_n}_{\gamma}:n\in\omega,\gamma\in U_{q_n}\}\] is a
constant sequence. So the top
level we are going to add does not affect those $b^p_{\gamma}$'s
for $\gamma\in U_q\smallsetminus U_{\omega}$.

Let $\{\langle x_m,\Gamma_m\rangle: m\in\omega\}$ be an enumeration
of $t'\times [U_{\omega}]^{<\omega}$. For each 
$\langle x_m,\Gamma_m\rangle$ we choose an increasing sequence
$\langle y_{m,i}:i\in\omega\rangle$ such that \[x_m=y_{m,0}<
y_{m,1}<\ldots,\] $y_{m,i+1}$ is a conservative extension
of $y_{m,i}$ with respect to $\Gamma_m$ and
\[\bigcup_{i\in\omega}ht(y_{m,i})=\alpha_q.\]
Now let $y_m=\bigcup_{i\in\omega}y_{m,i}$ and let $t_q=t'\cup\{y_m:
m\in\omega\}$. It is easy to see that $t_q\in {\cal T}$.
For each $\gamma\in U_{\omega}$ we define $b^q_{\gamma}$ to be an
extension of $b'_{\gamma}$ on $t_q$ such that
\[b^q_{\gamma}(y_m)= \bigcup_{i\in\omega}b'_{\gamma}(y_{m,i})\] 
for all $m\in\omega$. We define also $k_q$ to be an extension of $k'$
on $t_q$ such that for each $m\in\omega$, the tree $k_q(y_m)$ is in
$\cal T$, $ht(k_q(y_m))=\alpha_q+1$, $k_q(y_m)$ is an end-extension
of $\bigcup_{i\in\omega}k'(y_{m,i})$ and $b^q_{\gamma}(y_m)\in 
k_q(y_m)$ for all $\gamma\in U_{\omega}$. It is easy to see now
that the element $q$ is in $\Bbb P$ and is a lower bound of $p_n$'s
 and $q_n$'s.  \quad $\Box$

\begin{lemma}

The forcing notion $\Bbb P$ satisfies $\kappa$-c.c..

\end{lemma}

\noindent {\bf Proof:}\quad
Let $\{p_{\eta}:\eta\in\kappa\}\subseteq {\Bbb P}$. By a cardinality
argument and $\Delta$-system lemma there is an $S\subseteq\kappa$,
$|S|=\kappa$ and there is a triple $\langle\alpha_0,t_0,k_0\rangle$
such that for every $\eta\in S$ \[\langle\alpha_{p_{\eta}},t_{p_{\eta}},
k_{p_{\eta}}\rangle = \langle\alpha_0,t_0,k_0\rangle,\] and
$\{U_{p_{\eta}}:\eta\in S\}$ forms a $\Delta$-system with the root
$U_0$. Furthermore, we can assume that for each $\gamma\in U_0$,
\[b^{p_{\eta}}_{\gamma}=b^{p_{\eta'}}_{\gamma}\mbox{ and } 
f^{p_{\eta}}_{\gamma}=f^{p_{\eta'}}_{\gamma}\] for any $\eta,\eta'\in S$.
Since there are at most
$(|\omega_1^{\leqslant\alpha_0}|^{|t_0|})^{\omega_1} = 2^{\omega_1}$
sequences of length $\omega_1$ of the functions from $t_0$ to 
$\omega_1^{\leqslant\alpha_0}$, there are $\eta,\eta'\in S$ such that
\[\{b^{p_{\eta}}_{\gamma}:\gamma\in U_{p_{\eta}}\smallsetminus U_0\}
\mbox{ and }
\{b^{p_{\eta'}}_{\gamma}:\gamma\in U_{p_{\eta'}}\smallsetminus 
U_0\}\] are same set of functions. It is easy to see now that the element
\[p=\langle\alpha_0,t_0,k_0,U_{p_{\eta}}\cup U_{p_{\eta'}},
B_{p_{\eta}}\cup B_{p_{\eta'}},F_{p_{\eta}}\cup F_{p_{\eta'}}\rangle\]
is a common lower bound of $p_{\eta}$ and $p_{\eta'}$. \quad $\Box$

\begin{lemma}

All cardinals between $\omega_1$ and $\kappa$ in $V$ are collapsed
in $V^{\Bbb P}$.

\end{lemma}

\noindent {\bf Proof:}\quad 
For any $\gamma\in\kappa$ let \[f_{\gamma}=\bigcup\{f^p_{\gamma}:p\in G
\mbox{ and }\gamma\in U_p\}\] where $G\subseteq {\Bbb P}$ is a 
$V$-generic filter. It is easy to check that $range(f_{\gamma})=\gamma$.
Also $dom(f_{\gamma})\subseteq\omega_1$. So in $V^{\Bbb P}$ we have
$|\gamma|\leqslant\omega_1$. \quad $\Box$

\bigskip

\noindent {\bf Remark:}\quad By Lemma 12, Lemma 13 and Lemma 14 we have
\[V^{\Bbb P}\models (2^{\omega}=\omega_1^{V}=\omega_1\mbox{ and }
2^{\omega_1}=\kappa =\omega_2).\]

\begin{lemma}

Let $G\subseteq {\Bbb P}$ be a $V$-generic filter and 
let $T_G=\bigcup\{t_p:p\in G\}$. 
Then $T_G$ is an $\omega$-distributive Aronszajn tree in $V[G]$.

\end{lemma}

\noindent {\bf Proof:}\quad
It is easy to see that $T_G$ is an $\omega_1$-tree. Suppose there is a
$p_0\in {\Bbb P}$ such that \[p_0\forces\dot{B}
\mbox{ is a branch of }T_G.\] We construct $p_0\geqslant q_0\geqslant
p_1\geqslant q_1\geqslant\ldots$ similar to the construction in Lemma 12
such that \[p_{n+1}\forces z_n\in\dot{B}\cap (t_{q_n})_{\alpha_{q_n}}\]
for some $z_n\in\omega_1^{\alpha_{q_n}}$. For constructing $q_{n+1}$
we use almost same method as in Lemma 12 except that we require
$q_{n+1}$ satisfy the following condition (g'):

\begin{quote}

For any $x\in t_{p_{n+1}}$ and $\Gamma\in [U_{n+1}]^{<\omega}$ 
(see Lemma 12 for the definition of $U_{n+1}$)
there are infinitely many $x'\in (t_{q_{n+1}})_{\alpha_{q_{n+1}}}$ 
such that $x'$ is a conservative extension of $x$ with respect to
$\Gamma$.

\end{quote}

\noindent This can be done just by stretching 
$t_{q_{n+1}}$ a little bit higher
and manipulating those $b^{q_{n+1}}_{\gamma}\!\res\!(t_{q_{n+1}}\smallsetminus
t_{p_{n+1}})$ for $\gamma\in U_{n+1}$ more carefully. Let $q$ be a lower
bound of $\langle q_n:n\in\omega\rangle$ constructed same as in Lemma 12
except that for any $\langle x_m,\Gamma_m\rangle$ the sequence
$\langle y_{m,i}:i\in\omega\rangle$ is chosen such that 
$\bigcup_{i\in\omega}y_{m,i}$ is different from 
$\bigcup_{n\in\omega}z_n$. This is guaranteed by the condition (g').
Now \[\bigcup_{n\in\omega}z_n\not\in (t_q)_{\alpha_q}.\] Hence
\[q\forces\dot{B}\subseteq t_q.\] This contradicts that $B$ is a branch
of $T_G$ in $V[G]$.

Next we prove that $T_G$ is $\omega$-distributive. Let 
${\Bbb Q}=\langle T_G,\leqslant'\rangle$ be the forcing notion by
reversing tree order ($\leqslant'\;=\;\geqslant_{T_G}$). Given any
$\tau\in 2^{\omega}$ in $V^{{\Bbb P}*{\dot{\Bbb Q}}}$. 
It suffices to show that $\tau\in V$.
We construct a decreasing sequence 
\[\langle p_0,\dot{x}_0\rangle\geqslant
\langle q_0,\dot{x}_0\rangle\geqslant
\langle p_1,\dot{x}_1\rangle\geqslant
\langle q_1,\dot{x}_1\rangle\geqslant\ldots\] in ${\Bbb P}*\dot{\Bbb Q}$
such that \[\langle p_0,\dot{x}_0\rangle\forces\dot{\tau}\mbox{ is
a function from }\omega\mbox{ to }2,\]
\[p_n\forces\dot{x}_n\in\omega_1^{\alpha_{p_n}},\]
\[q_n\forces\dot{\tau}(n) = l_n\] for some $l_n\in\{0,1\}$ and
\[q_n\forces\dot{x}_n =\bar{x}_n\] for some $\bar{x}_n\in 
(t_{p_n})_{\alpha_{p_n}}$. In addition we can extend $q_n$ so that
the requirements for Player $I\!I$ to win the game are also satisfied.
Now we can construct a lower bound $q$ of $q_n$ same as we did in Lemma 12
except that we put also $x=\bigcup_{n\in\omega}\bar{x}_n$ into the top
level of $t_q$. It is easy to see that 
$\langle q,x\rangle\in {\Bbb P}*\dot{\Bbb Q}$ and there is a $\sigma=
\langle l_0,l_1,\ldots\rangle\in 2^{\omega}$ in $V$ such that 
\[\langle q,x\rangle\forces\dot{\tau} =\sigma.\quad \Box\]

\begin{lemma}

Let $G\subseteq {\Bbb P}$ be a $V$-generic filter and let 
$k_G=\bigcup\{k_p:p\in G\}$. Let $T_G$ and $\Bbb Q$ be same 
as in Lemma 15. Suppose $H\subseteq {\Bbb Q}$ is a 
$V[G]$-generic filter. 
Then $K_H=\bigcup\{k_G(x):x\in H\}$ is a Kurepa tree
in $V[G][H]$.

\end{lemma}

\noindent {\bf Proof:}\quad
It is easy to see that $K_H$ is an $\omega_1$-tree. For any 
$\gamma\in\kappa$ let \[b_{\gamma}=\bigcup\{b^p_{\gamma}:
p\in G\mbox{ and }\gamma\in U_p\}.\] Then $b_{\gamma}$ is a function
with domain $T_G$. Let \[W_{\gamma}=\bigcup\{b_{\gamma}(x):x\in H\}.\]
Then it is easy to see that $W_{\gamma}$ is a branch of $K_H$.
We need now only to show that $W_{\gamma}$ and $W_{\gamma'}$ are
different branches for different $\gamma,\gamma'\in\kappa$.
Given distinct $\gamma$ and $\gamma'$ in $\kappa$. Let
\[
\begin{array}{ll}

D^5_{\gamma,\gamma'}= & \{p\in {\Bbb P}:\beta^p_{\gamma}=
\beta^p_{\gamma'}=\alpha_p\mbox{ and }\\
 &(\forall x\in 
t_p\!\res\!\alpha_p)(\exists y\in t_p)(y\geqslant x\mbox{ and }
b^p_{\gamma}(y)\not=b^p_{\gamma'}(y))\}.

\end{array}
\]

{\bf Claim 16.1}\quad The set $D^5_{\gamma,\gamma'}$ is dense
in $\Bbb P$.

Proof of Claim 16.1:\quad
Given $p_0\in {\Bbb P}$. Without loss of generality we assume that 
$p_0\in D^2_{\gamma}\cap D^2_{\gamma'}$ and $\beta^{p_0}_{\gamma}
=\beta^{p_0}_{\gamma'}=\alpha_{p_0}$. First, we extend $t_{p_0}$
to $t_p\in {\cal T}$ such that 
\[\alpha_p=ht(t_p)=\alpha_{p_0}+\omega+1.\]
Then, we choose one extension $k_p$ of $k_{p_0}$ on $t_p$. 
Now we can easily extend 
$b^{p_0}_{\gamma}$ and $b^{p_0}_{\gamma'}$ to $b^p_{\gamma}$ and 
$b^p_{\gamma'}$ on $t_p$ while keeping other things unchanged
such that the resulting element $p$ is in
$\Bbb P$ and for each $x\in t_p\res\alpha_p$ there is an $y\in
(t_p)_{\alpha_p}$ and $y>x$ such that
$b^p_{\gamma}(y)\not=b^p_{\gamma'}$. It is easy to see the element
 $p$ is less than $p_0$ and is in $D^5_{\gamma,\gamma'}$. This ends the
proof of the claim.

\medskip

We need to prove $W_{\gamma}$ and $W_{\gamma'}$ are different branches
of $K_H$ in $V[G][H]$. Suppose $x\in H$ and \[x\forces\dot{W}_{\gamma}=
\dot{W}_{\gamma'}\] in $V[G]$. Let $p_0\in G$ be such that 
$x\in t_{p_0}$. By the claim we can find a $p\leqslant p_0$ and
$p\in G\cap D^5_{\gamma,\gamma'}$ such that $\alpha_p>ht(x)$.
Then we can choose $y\in t_p$ and $y>x$ such that $b^p_{\gamma}(y)\not=
b^p_{\gamma'}(y)$. Therefore 
\[y\forces \dot{W}_{\gamma}\not=\dot{W}_{\gamma'},\]
which contradicts that \[x\forces \dot{W}_{\gamma}=\dot{W}_{\gamma'}.
\quad \Box\]

\bigskip

The next lemma is probably the hardest part of this section.

\begin{lemma}

There are no Kurepa trees in $V^{\Bbb P}$.

\end{lemma}

\noindent {\bf Proof:}\quad
Suppose \[\forces_{\Bbb P}\dot{T}\mbox{ is a Kurepa tree with }\kappa
\mbox{ branches }\dot{\cal C}=\{\dot{c}_{\gamma}:\gamma\in\kappa\}.\]
For each $\gamma\in\kappa$ such that $cof(\gamma)=(2^{\omega_1})^+$
we choose an elementary submodel ${\goth A}_{\gamma}$ of $H(\lambda)$
such that 

(a) $|{\goth A}_{\gamma}|\leqslant 2^{\omega_1}$,

(b) $\{\dot{T},\dot{\cal C},{\Bbb P},\gamma\}\subseteq {\goth A}_{\gamma}$,

(c) $[{\goth A}_{\gamma}]^{\leqslant\omega_1}\subseteq {\goth A}_{\gamma}$.

\noindent By the Pressing Down Lemma we can find a set \[S\subseteq
\{\gamma\in\kappa:cof(\gamma)=(2^{\omega_1})^+\}\] with $|S|=\kappa$
such that 

(d) $\{{\goth A}_{\gamma}:\gamma\in S\}$ forms a $\Delta$-system with
the common root $\goth B$,

(e) there is a $\eta_0\in\kappa$ such that
$\eta_0=\bigcup\{\eta\in\kappa:\eta\in {\goth A}_{\gamma}\cap\gamma\}$
for every $\gamma\in S$,

(f) for any $\gamma,\gamma'\in S$ there is an isomorphism
$h_{\gamma,\gamma'}$ from ${\goth A}_{\gamma}$ to ${\goth A}_{\gamma'}$
such that $h_{\gamma,\gamma'}\res {\goth B}$ is an identity map.

Notice that $\omega_1\subseteq {\goth B}$ and $\omega_1^{<\omega_1}
\subseteq {\goth B}$. So for any $x\in\omega_1^{<\omega_1}$ we have
$h_{\gamma,\gamma'}(x)=x$. Let $\gamma_0$ be the minimal ordinal in $S$.
For any $p,p'\in {\Bbb P}$ we write $p\res {\goth A}_{\gamma} = p'$
to mean $\langle\alpha_p,t_p,k_p\rangle =
\langle\alpha_{p'},t_{p'},k_{p'}\rangle$, $U_p\cap {\goth A}_{\gamma}=
U_{p'}$, $b^p_{\gamma}=b^{p'}_{\gamma}$ and $f^p_{\gamma}=f^{p'}_{\gamma}$  
for each $\gamma\in U_{p'}$. 
We write also $p\res {\goth B}=p'$ to mean the same thing as above except
replacing ${\goth A}_{\gamma}$ by $\goth B$. Notice that
for $p,p'\in {\goth A}_{\gamma}$ the sentence $p\res {\goth B} = p'$
is first-order with parameters in ${\goth A}_{\gamma}$, {\em i.e.}
the term $\goth B$ could be eliminated. Next we are going to do
a complicated inductive construction of several sequences.

We construct inductively the sequences 

$\langle p_n\in {\Bbb P}:n\in\omega\rangle$,

$\langle p_s\in {\Bbb P}:s\in2^{<\omega}\rangle$,

$\langle \eta_n\in\omega_1:n\in\omega\rangle$ and

$\langle x_s\in\omega_1^{<\omega_1}:s\in2^{<\omega}\rangle$

\noindent in ${\goth A}_{\gamma_0}$ such that

(1) $p_{n+1}<p_n$ and $\alpha_{p_n}<\alpha_{p_{n+1}}$ for every
$n\in\omega$,

(2) $p_s\leqslant p_{s'}$ for any $s,s'\in 2^{<\omega}$ and 
$s'\subseteq s$,

(3) $p_s\res {\goth B}=p_n$ for any $n\in\omega$ and $s\in 2^n$,

(4) $\eta_n <\eta_{n+1}$ for every $n\in\omega$,

(5) $x_{s'}\leqslant x_s$ for any $s,s'\in 2^{<\omega}$ and 
$s'\subseteq s$,

(6) $ht(x_s)=\eta_n$ for any $s\in 2^n$,

(7) $x_s\not= x_{s'}$ for any $s,s'\in 2^n$ and $s\not= s'$,

(8) $p_s\forces x_s\in\dot{c}_{\gamma_0}$ for every $s\in 2^{<\omega}$,

(9) $t_{p_n}$ is consistent with respect to $\{b^{p_s}_{\gamma}:
\gamma\in\bigcup_{s\in 2^n}U_{p_s}\}$ for each $n\in\omega$,

(10) $\beta^{p_s}_{\gamma}=\alpha_{p_s}$ for all $\gamma\in U_{p_s}$
such that $\beta^{p_{s'}}_{\gamma}\not=\beta^{p_{s''}}_{\gamma}$ for 
some $s'\subseteq s''\subseteq s$,

(11) $\{b^{p_s}_{\gamma}:\gamma\in U_{p_s}\smallsetminus U_{p_n}\}$ and
$\{b^{p_{s'}}_{\gamma}:\gamma\in U_{p_{s'}}\smallsetminus U_{p_n}\}$
are the same set of functions for all $s,s'\in 2^n$.

\medskip

We need to add more requirements for those sequences along the 
inductive construction. 

For any $s\in 2^{<\omega}$ let \[U_s =\{\gamma\in U_{p_s}:\exists
s',s'' (s'\subseteq s''\subseteq s\mbox{ and }\beta^{p_{s'}}_{\gamma}
\not= \beta^{p_{s''}}_{\gamma})\}.\] Let's fix an onto function
$j:\omega\mapsto\omega\times\omega$ such that $j(n)=\langle a,b\rangle$
implies $a\leqslant n$. Let $\pi_1,\pi_2$ be projections from
$\omega\times\omega$ to $\omega$ such that $\pi_1(\langle a,b\rangle)=a$
and $\pi_2(\langle a,b\rangle)=b$. Let 
\[\xi_n:\omega\mapsto t_{p_n}\times ([\bigcup_{s\in 2^n}U_s]^{<\omega})\] 
and \[\zeta_n:\omega\mapsto\bigcup_{s\in 2^n}U_s\] 
be two onto functions for
each $n\in\omega$. Let $e$ be a function with $domain(e) = \omega$
such that \[e(n)=\xi_{\pi_1(j(n))}(\pi_2(j(n))).\]
The functions $\xi_n$'s, $\zeta_n$'s and $e$ are going to be used for
bookkeeping purpose. 
For $s\in 2^m$ and $m<n$ let \[C_{s,n}=\{s'\in 2^n:s\subseteq s'\}.\]
For any $m,n\in\omega$, $m\leqslant n$ let

\[
\begin{array}{ll}

Z^n_m = & \{b^{p_{s'}}_{\gamma}:s\in 2^{\pi_1(j(m))},\,\gamma\in
\pi_2(e(m))\cap U_s\mbox{ and }s'\in C_{s,n}\}\cup\\
  &\{b^{p_{s'}}_{\gamma}: s\in 2^{\pi_1(j(m))},\,\gamma\in U_s\mbox{ and }
\gamma=\zeta_{\pi_1(j(m))}(i)\mbox{ for some }i\leqslant n\}.

\end{array}
\]
Note that $Z^n_m$ is finite and for each
$b^p_{\gamma}\in Z^n_m$ we have $\beta^p_{\gamma}=\alpha_{p_n}$.
For each $m,n\in\omega$ we need also construct another set
\[Y^n_m =\{y_{m,i}:m\leqslant i\leqslant n\}.\] Then 
$Z^n_m$'s and $Y^n_m$'s and other four sequences should
satisfy two more conditions.

(12) $y_{m,m}=\pi_1(e(m))$ and $y_{m,i}\in (t_{p_i})_{\alpha_{p_i}}$
for $m<i\leqslant n$,

(13) $y_{m,i+1}$ is a conservative extension of $y_{m,i}$ with
respact to $Z^{i+1}_m$.

\noindent Next we do the inductive construction.
Suppose we have had sequences

$\langle p_n\in {\Bbb P}:n<l\rangle$,

$\langle p_s\in {\Bbb P}:s\in2^{<l}\rangle$,

$\langle \eta_n\in\omega_1:n<l\rangle$,

$\langle x_s\in\omega_1^{<\omega_1}:s\in 2^{<l}\rangle$,

$\{Z^n_m:n<l,m\leqslant n\}$ and 

$\{Y^n_m:n<l,m\leqslant n\}$.

\medskip

We first choose distinct $\{\gamma_s:s\in 2^l\}\subseteq S$.
For any $s\in 2^l$ let $p^s =h_{\gamma_0,\gamma_s}(p_{s\res l})$.
Note that \[p^s=\langle\alpha_{p_{s\res l}},t_{p_{s\res l}},
k_{p_{s\res l}}, U_{p^s},B_{p^s},F_{p^s}\rangle\]
where \[U_{p^s}=\{h_{\gamma_0,\gamma_s}(\gamma):\gamma\in 
U_{p_{s\res l}}\},\]
\[B_{p^s}=\{b^{p^s}_{h_{\gamma_0,\gamma_s}(\gamma)}:
b^{p^s}_{h_{\gamma_0,\gamma_s}(\gamma)}=
h_{\gamma_0,\gamma_s}(b^{p_{s\res l}}_{\gamma})\mbox{ and }
\gamma\in U_{p_{s\res l}}\}\] and \[F_{p^s}=\{h_{\gamma_0,\gamma_s}
(f^{p_{s\res l}}_{\gamma}):\gamma\in U_{p_{s\res l}}\}.\] Notice that
$b^{p^s}_{h_{\gamma_0,\gamma_s} (\gamma)}$ and 
$b^{p_{s\res l}}_{\gamma}$ are same functions with different indices.
Notice also that \[\alpha_{p_{s\res l}}=\alpha_{p_{l-1}},\,
t_{p_{s\res l}}=t_{p_{l-1}},\, k_{p_{s\res l}}=k_{p_{l-1}}\] and
\[U_{p^s}=\{h_{\gamma_0,\gamma_s}(\gamma):\gamma\in U_{p_{s\res l}}
\smallsetminus U_{p_{l-1}}\}\cup U_{p_{l-1}}.\]
Let \[\bar{p}_{l-1}=\langle\alpha_{p_{l-1}},t_{p_{l-1}},k_{p_{l-1}},
U_{\bar{p}_{l-1}},B_{\bar{p}_{l-1}},F_{\bar{p}_{l-1}}\rangle\] where 
\[U_{\bar{p}_{l-1}}=\bigcup_{s\in 2^l} U_{p^s},\]
\[B_{\bar{p}_{l-1}}=\{ b^{p^s}_{\gamma}:s\in 2^l,\gamma\in U_{p^s}\}\]
and \[F_{\bar{p}_{l-1}}=\{f^{p^s}_{\gamma}:s\in 2^l, \gamma\in
U_{p^s}\}.\] Since $t_{p-1}$ is consistent with $\bigcup_{s\in 2^{l-1}}
B_{p_s}$ by (9), then we have $\bar{p}_{l-1}\in {\Bbb P}$.
Since \[\bar{p}_{l-1}\forces\{\dot{c}_{\gamma_s}:s\in 2^l\}\mbox{ is a
set of distinct branches of }\dot{T},\] then there exist $\bar{p}_l
\leqslant \bar{p}_{l-1}$, $\eta_l\in\omega_1$ such that
$\eta_l>\eta_{l-1}$, and there exist distinct 
\[\{x_s:s\in 2^l\}\subseteq \omega^{\eta_l}\]
 such that \[\bar{p}_l\forces x_s\in\dot{c}_{\gamma_s}\]
for all $s\in 2^l$. We can also require that $\alpha_{\bar{p}_l}>
\alpha_{\bar{p}_{l-1}}$ and $\beta^{\bar{p}_l}_{\gamma}=
\alpha_{\bar{p}_l}$ for all $\gamma\in U_{\bar{p}_l}$ such that
$\beta^{\bar{p}_l}_{\gamma}>\beta^{\bar{p}_{l-1}}_{\gamma}$,
or for all $\gamma\in
\bigcup_{s\in 2^{l-1}}h_{\gamma_0,\gamma_s}[U_{s\res l}]$.
For each $s\in 2^l$ let $\bar{U}_s$ be a set of
$\omega_1$ ordinals such that $\bar{U}_s\subseteq
{\goth A}_{\gamma_s}\smallsetminus {\goth B}$ and $\bar{U}_s\cap 
U_{\bar{p}_l}=\emptyset$. 
Since $B_{\bar{p}_l}$ has only $\leqslant\omega_1$ functions, 
we can use the ordinals in $\bar{U}_s$ to re-index all functions
in $B_{\bar{p}_l}$, say $B_{\bar{p}_l}$ and $\{b^{\bar{p}_l}_{\gamma}:
\gamma\in\bar{U}_s\}$ are same set of functions. 
Let $f^{\bar{p}_l}_{\gamma}$ be an empty
function for each $\gamma\in \bar{U}_s$. We now construct a $\bar{p}$
such that 
\[\bar{p}=\langle\alpha_{\bar{p}_l},t_{\bar{p}_l},k_{\bar{p}_l},
U_{\bar{p}},B_{\bar{p}},F_{\bar{p}}\rangle,\] where
\[U_{\bar{p}}=U_{\bar{p}_l}\cup (\bigcup_{s\in 2^l}\bar{U}_s),\]
\[B_{\bar{p}}=B_{\bar{p}_l}\cup (\bigcup_{s\in 2^l}
\{b^{\bar{p}_l}_{\gamma}:\gamma\in\bar{U}_s\})\] and
\[F_{\bar{p}}=F_{\bar{p}_l}\cup (\bigcup_{s\in 2^l}
\{f^{\bar{p}_l}_{\gamma}:\gamma\in\bar{U}_s\}).\] 
It is easy to see that $\bar{p}\in {\Bbb P}$ and $\bar{p}\leqslant
\bar{p}_l$.

\bigskip

{\bf Claim 16.2}\quad For each $s\in 2^l$ let $\bar{p}_s=
\bar{p}\res {\goth A}_{\gamma_s}$. Then 
$\bar{p}_s\forces x_s\in\dot{c}_{\gamma_s}$.

Proof of Claim 16.2:\quad
It is true that $\bar{p}_s\in {\goth A}_{\gamma_s}$ because
$({\goth A}_{\gamma_s})^{\leqslant\omega_1}\subseteq 
{\goth A}_{\gamma_s}$.
Suppose \[\bar{p}_s\not\forces x_s\in\dot{c}_{\gamma_s}.\] Then there is a
$p'_s\leqslant\bar{p}_s$ such that \[p'_s\forces x_s\not\in\dot{c}_{\gamma_s}.\]
Since ${\goth A}_{\gamma_s}\preceq H(\lambda)$, we can choose $p'_s\in
{\goth A}_{\gamma_s}$. It is now easy to see that $p'_s$ and $\bar{p}$
are compatible (here we use the fact that every function in
$B_{\bar{p}}$ is also in $B_{\bar{p}_s}$ with possibly different index).
This derives a contradiction.

\medskip

Let $p_l=\bar{p}\res {\goth B}$ and $p_s=
h_{\gamma_s,\gamma_0}(\bar{p}_s)$. Then

$\langle p_n:n\leqslant l\rangle$,

$\langle p_s:s\in 2^{\leqslant l}\rangle$,

$\langle \eta_n:n\leqslant l\rangle$ and

$\langle x_s:s\in 2^{\leqslant l}\rangle$ 

\noindent satisfy conditions
(1)---(11). For example, we have \[p_s\forces x_s\in\dot{c}_{\gamma_0}\]
because $p_s =h_{\gamma_s,\gamma_0}(\bar{p}_s)$, $\gamma_0=
h_{\gamma_s,\gamma_0}(\gamma_s)$ and 
$x_s=h_{\gamma_s,\gamma_0}(x_s)$. We have also that $t_{p_l}$ is
consistent with $\{b^{p_s}_{\gamma}:
s\in 2^l\mbox{ and }\gamma\in U_{p_s}\}$ because $\bar{p}\in {\Bbb P}$.

We need to deal with the conditions (12) and (13).

For each $m\leqslant l$ the set $Z^l_m$ has been defined before.
For $m<l$ since $Z^l_m$ is finite, there exists a $y_{m,l}\in 
(t_{p_l})_{\alpha_{p_l}}$ such that $y_{m,l}$ is a consistent extension
of $y_{m,l-1}$ with respect to $Z^l_m$. Let $y_{l,l}=\pi_1(e(l))$.
It is not hard to see that those sequences up to stage $l$ satisfy
conditions (12) and (13). This ends the construction.

\medskip

We want to draw the conclusion now.

For each $m\in\omega$ let $y_m=\bigcup_{i\in\omega}y_{m,i}$ and let
\[t_{p_{\omega}}=(\bigcup_{n\in\omega}t_{p_n})\cup\{y_m:m\in\omega\}.\]
It is easy to see that $t_{p_{\omega}}\in {\cal T}$. Let
$\alpha_{p_{\omega}}=\bigcup_{n\in\omega}\alpha_{p_n}$. Then
$ht(t_{p_{\omega}})=\alpha_{p_{\omega}}+1$. Let 
\[U=\bigcup_{s\in 2^{<\omega}}U_s=\{\gamma:\exists\tau\in 2^{\omega}
\mbox{ such that }\bigcup\{\beta^{p_{\tau\res n}}_{\gamma}:n\in\omega
\mbox{ and }
\gamma\in U_{p_{\tau\res n}}\}=\alpha_{p_{\omega}}\}.\]
Then $U$ is a countable set. Notice that for any 
$s\in 2^{<\omega}$ and $\gamma\in U_{p_s}\smallsetminus U$,
for any $s'\supseteq s$ we have $\beta^s_{\gamma}=\beta^{s'}_{\gamma}$.
Let $k'=\bigcup_{n\in\omega}k_{p_n}$. For each $\tau\in 2^{\omega}$ and
$\gamma\in\bigcup_{s\in 2^{<\omega}}U_{p_s}$ let
\[b^{\tau}_{\gamma}=\bigcup\{b^{p_{\tau\res n}}_{\gamma}:\gamma\in
U_{p_{\tau\res n}},n\in\omega\}\] and let
\[f^{\tau}_{\gamma}=\bigcup\{f^{p_{\tau\res n}}_{\gamma}:\gamma\in
U_{p_{\tau\res n}},n\in\omega\}.\]
For each $m\in\omega$ and $\gamma\in U$ we define
\[b^{\tau}_{\gamma}(y_m)=\bigcup\{b^{\tau}_{\gamma}(y_{m,i}):
m\leqslant i<\omega\}.\]
Since for each $\gamma\in U$ and $m\in\omega$ 
there exists an $n$ such that for any
$s,s'\in 2^l$ for $l\geqslant n$ and $s\!\res\! n=s'\!\res\! n$ we have
$b^{p_s}_{\gamma}(y_m)=b^{p_{s'}}_{\gamma}(y_m)$. 
This is guaranteed by the construction of $Z^n_m$'s and $Y^n_m$'s. So
for any $\tau,\tau'\in 2^{\omega}$, \[\tau\!\res\! n=\tau'\!\res\! n
\mbox{ implies }b^{\tau}_{\gamma}(y_m)=b^{\tau'}_{\gamma}(y_m).\]
 Hence for each $m\in\omega$ the set
\[\{b^{\tau}_{\gamma}(y_m):\tau\in 2^{\omega},\gamma\in U\}\] is
countable. (This is why the condition (g) of Definition 10 is needed.)
Let $k'(y_m)$ be in $\cal T$ such that \[\bigcup_{i\in\omega}
k'(y_{m,i})\leqslant_{end}k'(y_m)\] and \[\{b^{\tau}_{\gamma}(y_m):
\tau\in 2^{\omega},\gamma\in U\}\subseteq
(k'(y_m))_{\alpha_{p_{\omega}}}.\] Then let $k_{p_{\omega}}=k'$.
For each $\tau\in 2^{\omega}$ let $x_{\tau}=\bigcup_{n\in\omega}
x_{\tau\res n}$. Then $x_{\tau}\in\omega_1^{\alpha_{p_{\omega}}}$.
For each $\tau\in 2^{\omega}$ let $p_{\tau}$ be the lower bound of
$\{p_{\tau\res n}:n\in\omega\}$ constructed same as in Lemma 12.
Then we have $p_{\tau}\in {\Bbb P}$ and \[p_{\tau}\forces x_{\tau}
\in\dot{c}_{\gamma_0}.\]
Choose distinct ordinals $\{\gamma_{\tau}:\tau\in O\}\subseteq S$ for some
$O\subseteq 2^{\omega}$ and $|O|=\omega_1$.
Let $p^{\tau}=h_{\gamma_0,\gamma_{\tau}}(p_{\tau})$.
Then \[p^{\tau}\forces x_s\in\dot{c}_{\gamma_{\tau}}\] for any 
$\tau\in O$.
Let \[q=\langle\alpha_{p_{\omega}},t_{p_{\omega}},k_{p_{\omega}},
U_q,B_q,F_q\rangle,\] where \[U_q=\bigcup_{\tau\in O}
\{h_{\gamma_0,\gamma_{\tau}}(\gamma):\gamma\in U_{p_{\tau}}\},\]
\[B_q=\bigcup_{\tau\in O}
\{h_{\gamma_0,\gamma_{\tau}}(b^{\tau}_{\gamma}):\gamma\in
U_{p_{\tau}}\}\] and \[F_q=\bigcup_{\tau\in O}\{
h_{\gamma_0,\gamma_{\tau}}(f^{\tau}_{\gamma}):\gamma\in
U_{p_{\tau}}\}.\]

\bigskip

{\bf Claim 16.3}\quad The element $q$ is in $\Bbb P$ and $q\leqslant
p^{\tau}$ for all $\tau\in O$.

Proof of Claim 16.3:\quad It is easy to see that
$|U_q|\leqslant\omega_1$ (the condition $|U_p|\leqslant\omega_1$ for 
$p\in {\Bbb P}$ in Definition 10 is needed here since if we require
only $|U_p|<\omega_1$, then $q$ wouldn't be in $\Bbb P$). 
It is also easy to see that for each $\tau\in O$ we have 
$q\!\res\! {\goth A}_{\gamma_{\tau}}=p^{\tau}$. Hence it suffices to show
that $t_{p_{\omega}}$ is consistent with $B_q$. But this is guaranteed
by condition (9) and the construction of $y_m$'s.

\bigskip

{\bf Claim 16.4}\quad
$q\forces (\dot{T})_{\alpha_{p_{\omega}}}\mbox{ is uncountable.}$

Proof of Claim 16.4:\quad
This is because of the facts $x_{\tau}\not=x_{\tau'}$ 
for different $\tau,\tau'\in O$, $|O|=\omega_1$,
\[q\forces x_{\tau}\in\dot{c}_{\gamma_{\tau}}\]
and \[q\forces\dot{c}_{\gamma_{\tau}}\subseteq\dot{T}.\]

\bigskip

By above claim we have derived a contradiction that
\[\forces(\dot{T}\mbox{ is a Kurepa tree)}\] but
\[q\forces(\dot{T}\mbox{ is not a Kurepa tree}). \quad \Box\]

\section{Questions}

We would like to ask some questions.

\begin{question}

Suppose our ground model is the  L\'{e}vy model defined in the first
section. Can we find a proper forcing notion
such that the forcing extension will contain Kurepa trees?
If the answer is `no', then we would like to know if there are
any forcing notions of size $\leqslant\omega_1$ which 
preserve $\omega_1$ such that the generic extension contains
Kurepa trees?

\end{question}

\begin{question}

Suppose the answer of one of the questions above
is {\em Yes}. Is it true that given any model of \ch\, 
there always exists an $\omega_1$-preserving forcing
notion of size $\leqslant\omega_1$ such that forcing
with that notion creates Kurepa trees in the generic extension?

\end{question}

\begin{question}

Does there exist a model of \ch\, plus no Kurepa trees,
in which there is a {\em c.c.c.}-forcing notion of size
$\leqslant\omega_1$ such that forcing with that notion 
creates Kurepa trees in the generic extension?
If the answer is {\em Yes}, then we would like to ask
the same question with {\em c.c.c.} replaced by
one of some nicer chain conditions such as $\aleph_1$-caliber,
Property $K$, etc.

\end{question}

\bigskip

Department of Mathematics, University of Illinois

Urbana, IL 61801, USA

{\em e-mail: jin@@math.uiuc.edu}

\bigskip

Institute of Mathematics, The Hebrew University

Jerusalem, Israel

\bigskip

Department of Mathematics, Rutgers University

New Brunswick, NJ, 08903, USA

\bigskip

{\em Sorting:} The first address is the first author's and the last
two are the second author's. 

\end{document}